%Paper: dg-ga/9411008
%From: huebschm@mpim-bonn.mpg.de (Johannes Huebschmann)
%Date: Wed, 23 Nov 1994 10:51:18 --100

%&amstex
%This paper can be processed with AMSTeX 2.1.
\documentstyle{amsppt}
\magnification=1200
\hoffset=-0.5pc
\vsize=57.2truepc
\hsize=38.5truepc
\nologo
\spaceskip=.5em plus.25em minus.20em

\define\abramars{1}
\define\armgojen{2}
\define\armamonc{3}
\define\atibottw{4}
\define\bocosroy{5}
\define\ebinone{6}
\define\eellsone{7}
\define\goldmone{8}
 \define\howeone{9}
 \define\singula{10}
\define\singulat{11}
\define\topology{12}
  \define\direct{13}
\define\singulth{14}
 \define\locpois{15}
 \define\modusym{16}
\define\huebjeff{17}
\define\kobanomi{18}
\define\kostafou{19}
\define\lermonsj{20}
\define\marswein{21}
 \define\mathone{22}
\define\mittvial{23}
\define\narashed{24}
\define\naramntw{25}
\define\nararama{26}
\define\raghuboo{27}
\define\gwschwar{28}
\define\gwschwat{29}
\define\sjamlerm{30}
\define\weilone{31}
\define\weiltwo{32}
\define\weylbook{33}
\define\whitnone{34}
\define\whitnfou{35}
\noindent
dg-ga/9411008
\topmatter
\title
Smooth structures on certain moduli spaces\\
for bundles on a surface
\endtitle
\author Johannes Huebschmann{\dag}
\endauthor
\affil
Max Planck Institut f\"ur Mathematik
\\
Gottfried Claren-Str. 26
\\
D-53 225 BONN
\\
Huebschm\@mipm-bonn.mpg.de
\endaffil
\date{June 23, 1994}
\enddate

\endkeywords

\endsubjclass
\abstract{
Let $\Sigma$ be a closed  surface, $G$ a compact Lie group, not necessarily
connected, with Lie algebra $g$, $\xi \colon P \to \Sigma$ a principal
$G$-bundle, let $N(\xi)$ denote the moduli space  of central Yang-Mills
connections on $\xi$, with reference to suitably chosen additional data, and
let $\roman{Rep}_{\xi}(\Gamma,G)$ be the  space of representations of the
 universal central extension $\Gamma$ of the fundamental group of $\Sigma$ in
$G$ that corresponds to $\xi$. We construct smooth structures on $N(\xi)$ and
$\roman{Rep}_{\xi}(\Gamma,G)$, that is, algebras of continuous functions
 which restrict to smooth functions on the strata of certain associated
stratifications; by means of a detailed investigation of the derivative of the
 holonomy we show thereafter that, with reference to these smooth structures,
 the assignment to a smooth connection $A$ of its holonomies  with reference
to suitable closed paths yields a diffeomorphism from $N(\xi)$ onto
$\roman{Rep}_{\xi}(\Gamma,G)$;  moreover we show that the derivative of the
 latter at the non-singular points of $N(\xi)$ amounts to a certain twisted
 integration mapping relating a suitable de Rham theory with group cohomology
 with appropriate coefficients. Finally we examine the infinitesimal geometry
of these moduli spaces with reference to the smooth structures and, for
 illustration, we show that, on the moduli space of flat
$\roman{SU}(2)$-connections for a surface of genus two which, as a space, is
 just complex projective 3-space, our smooth structure looks rather different
 from the standard structure.}
\endabstract

\thanks{{\dag} The author carried out this work in the framework
of the VBAC research group of Europroj.}
\endthanks
\endtopmatter
\document
\rightheadtext{Smooth structures on moduli spaces}
\leftheadtext{Johannes Huebschmann}

\beginsection Introduction

Let $X$ be a
decomposed topological space,
each piece
of the decomposition
being a smooth manifold.
A {\it smooth structure\/}
on $X$
is an algebra $C^{\infty}(X)$
of continuous functions
on $X$
which,
on each piece,
restrict to smooth functions.
We shall refer to such a space
as a {\it smooth\/} space.
In the present paper, we endow certain moduli spaces
with a smooth structure and thereafter analyze
their
singular structure and
infinitesimal geometry
by means of it.
It belongs to a series of papers about  a program
revealing the structure of these moduli spaces
by means of the symplectic or more generally Poisson geometry
of certain related classical constrained systems
but its results are of interest
in their own right.
In \cite\direct\
we construct the searched for Poisson structures on the moduli spaces,
thereby
obtaining structures of a {\it stratified symplectic space\/}
in the sense of \cite\sjamlerm;
such a structure encapsulates
the {\it mutual positions\/} of symplectic structures
on  the strata.
It is known that
some of these moduli spaces
carry the additional structure
of a (complex) projective
variety
which, however,
does not shed too much light
on the singular behaviour of
the symplectic
or Poisson
structures in general;
in fact it may happen that the symplectic structure
is singular whereas the complex analytic one is not.
An example will be mentioned shortly.
On the other hand, the
singular behaviour of
the symplectic
or more generally Poisson structures
can entirely be
understood in the framework of the real algebraic geometry
of appropriate smooth structures on these moduli spaces,
to which the present paper is devoted.
We shall relate the smooth structures
with appropriate complex analytic structures
elsewhere by means of a suitable notion of polarization
for Poisson  structures;
this will generalize the classical description of
a K\"ahler structure in terms of a holomorphic polarization
and in particular will provide the necessary means to
talk about {\it mutual positions\/} of K\"ahler structures
on the strata.
\smallskip
We explain at first briefly the moduli spaces.
Let
$\Sigma$ be a closed surface,
$G$ a compact
Lie group,
not necessarily connected,
with Lie algebra $g$,
and $\xi \colon P \to \Sigma$
a principal $G$-bundle,
having a connected total space $P$.
Further, pick
a Riemannian metric on $\Sigma$ and
an {\it orthogonal structure\/}
on $g$,
that is, an adjoint action invariant scalar product.
These data then determine a Yang-Mills
theory
studied
for connected $G$
extensively by {\smc Atiyah-Bott}
in \cite\atibottw\ to which we refer for
background and notation.
We only mention that
a connection is said to be
{\it Yang-Mills\/} provided it satisfies
the corresponding
Yang-Mills equations
and {\it central\/}
when its curvature is a 2-form on $\Sigma$
with values in the Lie algebra of the centre of $G$.
The {\it moduli space\/} $N(\xi)$
of central Yang-Mills
connections
is
then that of gauge equivalence classes of
central Yang-Mills
connections;
it is a compact space,
including as special cases
certain moduli spaces of flat connections and
the {\smc Narasimhan-Seshadri}-moduli spaces \cite\narashed\
of semi-stable holomorphic vector bundles.
For example, as a space,
the moduli space of flat
$\roman{SU}(2)$-connections for a surface
of genus 2 is just complex projective 3-space \cite\naramntw;
as a complex analytic space, it is non-singular but the symplectic
or more general stratified symplectic structure
degenerates on a Kummer surface; see~\cite\locpois.
For a general bundle $\xi$ and structure group $G$, we shall assume throughout
that the space $N(\xi)$ is non-empty, that is, that Yang-Mills
connections exist. For example, this will be the case for a connected
structure group, cf.~\cite\atibottw.
\smallskip
In \cite\topology\  we  have shown that
the assignment to a connection of its holonomies,
with reference to suitably chosen closed paths,
induces a
homeomorphism,
referred to as {\it Wilson loop mapping\/}
for a reason given in Section 2 below,
from $N(\xi)$
onto  a certain
representation space
$\roman {Rep}_{\xi}(\Gamma,G)$
for the universal central extension
$\Gamma$ of the fundamental group $\pi$
of $\Sigma$.
While the space $N(\xi)$
{\it depends\/}
on the choice of Riemannian metric
on $\Sigma$
the space
$\roman{Rep}_{\xi}(\Gamma,G)$
does not.
One of our aims is to show
that, with reference to
appropriate additional structure,
the  Wilson loop mapping
is in fact
a diffeomorphism.
\smallskip
We now give a brief overview of the paper.
Section 1 is preliminary in character.
In Section 2 we determine the derivative
of the holonomy,
viewed as a map from
the space of connections
to the structure group, once
the appropriate
additional requisite data
have been chosen.
In Section 3 we
introduce our algebras
of smooth functions
and
spell out the {\it first chief result\/} of the paper, Theorem 3.8;
it will say that,
the spaces
$N(\xi)$ and
$\roman {Rep}_{\xi}(\Gamma,G)$
being decomposed into connected components of orbit types
in the appropriate sense, the  Wilson loop mapping
is fact
a diffeomorphism.
In Section 4 we
give a
description of  the twisted integration mapping
tailored to our purposes.
In Section 5 we
rework
and extend
the classical relationship
between the
infinitesimal structure
of
representation spaces
and group cohomology
which goes back at least to
{\smc Weil}~\cite\weilone,\ \cite\weiltwo,\
cf.
\cite\raghuboo.
In Section 6
we reduce the smooth structures
of $N(\xi)$
and
$\roman {Rep}_{\xi}(\Gamma,G)$
near any of its points
to that of local models
of a kind
introduced in an
earlier paper \cite\singula,
endowed with suitable smooth structures.
This will be our {\it second chief result\/}.
Our {\it third chief result\/},
(6.15) below,
will be the existence of
suitable
partitions of unity;
this will then enable
us to complete the proof of Theorem 3.8
mentioned above.
In Section 7 we
examine the infinitesimal structure of our spaces
of interest.
In particular,
we shall establish the fact that the space
$N(\xi)$ is locally semi-algebraic.
Finally in Section 8 we examine the
moduli space of flat $\roman{SU}(2)$-connections for
a surface of genus 2 which,
cf. what was said above,
as a space is just complex projective
3-space.
We shall see that, as a smooth space with the appropriate
smooth structure, it looks rather different;
for example, at 16 isolated points,
the Zariski tangent space
has (real) dimension 10.
\smallskip
Abstracting the structure
of the spaces
$N(\xi)$ and
$\roman {Rep}_{\xi}(\Gamma,G)$
isolated in the present paper
we are led to
spaces with an algebra
of functions
which, locally, look like the reduced space
of
a momentum mapping
for a representation of a compact Lie group which
varies over the space,
with the obvious smooth structure on the reduced space.
This class of spaces may well be worth
an independent investigation.
\smallskip

I am indebted to
A. Weinstein
for discussions at various stages
of the project.
In particular, the final versions
of (7.9) and (7.10) below were found
after discussions with him.

\medskip\noindent{\bf 1. Preliminaries}\smallskip\noindent
Let $M$ be
a finite dimensional
smooth compact connected manifold,
$G$ a
(real) Lie group, $g$ its (real) Lie algebra, and
$\xi \colon P \to M$
a principal $G$-bundle over $M$,
with $G$ acting on the right of $P$.
We denote the action of
$x \in G$ by
$\roman R_x\colon p \mapsto px$, where $p \in P$.
The affine space $\Cal A(\xi)$
of
smooth connections on $\xi$
inherits an obvious action of
the group $\Cal G(\xi)$
of gauge transformations
and so
does the graded vector space $\Omega^*(M,\roman{ad}(\xi))$.
We pick a   {\it base point\/}
$Q \in M$ and a pre-image
$\widehat Q \in P$;
then assignment
to a gauge transformation $\gamma$ on $\xi$
of
$x_{\gamma} \in G$
defined
by $\gamma(\widehat Q) = \widehat Q x_{\gamma}$
furnishes a surjective homomorphism
$$
\Cal G(\xi)
@>>>
G
\tag1.1
$$
whose kernel
is the {\it group\/}
$\Cal G^Q(\xi)$
of (at $Q$)
{\it based  gauge transformations\/}.
The {\it adjoint bundle\/}
$\roman{ad}(\xi)$
is the Lie algebra bundle
over $M$
associated with $\xi$ and the adjoint action of $G$ on $g$.
With the obvious bracket,
its space of sections
$\Omega^0(M,\roman{ad}(\xi))$
is the Lie algebra
of $\Cal G(\xi)$ in a natural fashion.
The tangent bundle of a smooth manifold $X$
will be written
$\tau_X\colon \roman TX \to X$.
\smallskip
We shall not distinguish in notation between the naive
objects and their Sobolev completions
\cite\atibottw,
\cite\ebinone,
\cite\mittvial,
\cite\nararama.

\medskip\noindent{\bf 2. The derivative of the holonomy}
\smallskip\noindent
Let
$I= [0,b]$ be an interval and
$u\colon I \to M$  a smooth
path in $M$ having
starting point $Q$.
For a
connection $A$, we denote by
${
u_{A,\widehat Q}\colon I @>>> P
}$
the {\it horizontal\/}
lift
of $u$,
having starting  point $\widehat Q$.
For $t \in I$, let
${
u_{A,\widehat Q,t}\colon [0,t] @>>> P
}$
be the restriction of
$u_{A,\widehat Q}$ to $[0,t]$.
\smallskip

Among the various descriptions of
the space
$\Omega^j(M,\roman{ad}(\xi))$
of $j$-forms with values in
the adjoint bundle
$\roman{ad}(\xi)$
we shall take here that in terms of $G$-invariant horizontal
$g$-valued forms on $P$.
The following will be crucial.

\proclaim{Theorem 2.1}
With reference to a suitable Sobolev topology
on $\Cal A(\xi)$,
the assignment to
$(A,t) \in \Cal A(\xi) \times I$
of the
horizontal lift
$u_{A,\widehat Q}(t)$
furnishes a continuous map
$U$ from
$\Cal A(\xi) \times I$ to $P$
whose restriction to any smooth finite
dimensional
submanifold of
$\Cal A(\xi) \times I$
is smooth.
Given
a
connection $A$ on $\xi$ and a 1-form
$\vartheta \in \Omega^1(M,\roman{ad}(\xi))=\roman T_A\Cal A(\xi)$,
an explicit formula for the partial derivative
${
\frac{\partial U}{\partial \vartheta}(A,t)
=
dU(A,t)(\vartheta,0)
}$
is given by
$$
\frac{\partial U}{\partial \vartheta}(A,t)
=
u_{A,\widehat Q}(t) \int_{u_{A,\widehat Q,t}} \vartheta
\in
\roman T_{U(A,t)}P.
\tag2.2
$$
\endproclaim

\smallskip

\noindent
{\smc Remark 2.3.}
Some comment about the interpretation of the formula (2.2) might be in order:
The 1-form
$\vartheta$ being viewed as a $G$-invariant
$g$-valued
one on $P$
which vanishes on the vertical vectors,
the integral
$\int_{u_{A,\widehat Q,t}} \vartheta$
is well defined as an element of the Lie algebra $g$.
Moreover, by construction,
$u_{A,\widehat Q}(t) \in P$, and
the expression
$$
u_{A,\widehat Q}(t) \int_{u_{A,\widehat Q,t}} \vartheta
\in
\roman T_{u_{A,\widehat Q}(t)}P
$$
refers to the element which is obtained when the
canonical injection
from
$P \times g$  into the total space $\roman T P$
is applied to the pair
$\left(u_{A,\widehat Q}(t), \int_{u_{A,\widehat Q,t}} \vartheta\right)$.
\smallskip

\noindent
{\smc Remark 2.4.}
The existence of the derivative
of $U$, restricted to an arbitrary
smooth finite dimensional
submanifold,
and that of corresponding derivatives of arbitrarily
high order, follows from standard
facts about  analytical dependence
of the solution of a differential equation
on suitable parameters.

\demo{Proof of {\rm (2.1)}}
In the good range
$k > \frac {\dim M}2$,
convergence in the Sobolev topology
$H^k$
implies uniform
convergence,
cf.
\cite\eellsone\ Section 6.
This implies readily that
$U$ is continuous.
\smallskip

The smooth tangent space
$\roman T_{A}\Cal A(\xi)$ is
naturally identified with the
vector space
of
1-forms
$\Omega^1(M,\roman{ad}(\xi))$
with values in the adjoint bundle,
and, for a fixed value of $t \in I$,
we look for the derivative
${
\roman T_{A} U_t \colon
\roman T_{A}\Cal A(\xi)
@>>>
\roman T_{U_t}P
}$
of the map
$
U_t \colon \Cal A(\xi) @>>> P
$
which is given by the assignment to
a connection $A$ of  the value $U_t(A)=u_{A,\widehat Q}(t) \in P$.
Thus,
given $\vartheta \in
\Omega^1(M,\roman{ad}(\xi))$,
all we need is an expression
for the partial derivative
$$
\frac{\partial U_t}{\partial \vartheta}(A)
= \roman T_{A}U_t(\vartheta)
\in
\roman T_{u_{A,\widehat Q}(t)}P.
$$
To obtain such an expression,
given $s \in \bold R$ and
$\vartheta \in
\Omega^1(M,\roman{ad}(\xi))$,
we consider the horizontal lift
${
u_{A +  s \vartheta, \widehat Q}
\colon I @>>> P
}$
of $u$.
It is clear that
the assignment
to $(s,t)\in I \times I$
of $u_{A + s \vartheta, \widehat Q}(t)$
yields a smooth
map
${
\widehat u \colon I \times I @>>> P,
}$
and what we are looking for is an expression for the
partial derivative of this map
at $s=0$, whatever $t \in I$.
To simplify notation,
write
$
v = u_{A , \widehat Q}
\colon I @>>> P
$
for
the horizontal lift of $u$.
It is obvious that there is a unique map
${
a \colon I \times I @>>> G
}$
such that, for every $(s,t) \in I \times I$,
$$
\widehat u(s,t) =
v(t) a(s,t).
$$
When we fix $s$ and
differentiate this identity with respect to the parameter $t$
we obtain
the identity
$$
\widehat u_t' = v_t'a_t + v_t a_t';
$$
here we have written
$a_t= a(s,t) \in G$,
$\widehat u_t = \widehat u(s,t) \in P$,
$v_t = v(t) \in P$;
furthermore,
with a notation used e.~g. on  p. 69 of
\cite\kobanomi,
$\widehat u_t'$ is the tangent vector to the curve
$(s,t) \mapsto \widehat u(s,t) $ ($s$ fixed)
at the point $u(s,t)$,
and
$v_t'$ and
$a_t'$
refer to the corresponding tangent vectors
of the other curves coming into play.
Let $\omega \colon \roman TP \to g$
be the connection form of $A$;
then
$\omega - s \vartheta$
is the connection form of
$A + s \vartheta$.
Exploitation of the fact that
$\widehat u_t'$
is horizontal
for
the connection
$A + s \vartheta$
yields
$$
\align
0 &=\left(\omega - s \vartheta\right) \left(\widehat u_t'\right)
\\
&=\omega\left(u_t'\right) - s \vartheta\left(\widehat u_t'\right)
\\
&=\omega\left(v_t'a_t + v_t a_t'\right) -
s \vartheta\left(v_t'a_t + v_t a_t'\right)
\\
&=\omega\left(v_t'a_t\right)
 + \omega\left(v_t a_t'\right) -
s \vartheta\left(v_t'a_t\right)
 -
s \vartheta\left(v_t a_t'\right)
\\
&=\omega\left(\left(\roman R_{a_t}\right)_*v_t'\right)
 + \omega\left(v_t a_t'\right) -
s \vartheta\left(\left(\roman R_{a_t}\right)_* v_t'\right),
\endalign
$$
since
$v_t a_t'$ is vertical and since
$\vartheta$ is zero on vertical vectors;
we remind the reader that
$R_{a_t}\colon P \to P$
refers to the action of $G$ on P.
Moreover,
since the curve $v_t$ is horizontal with respect to $A$,
$G$-invariance of
$\omega$
implies that
$
\omega\left(\left(\roman R_{a_t}\right)_*v_t'\right)
$
equals
$\roman{ad}_{a_t^{-1}}
\omega\left(v_t'\right)$
which is zero;
likewise,
$G$-invariance of
$\vartheta$
implies
${
\vartheta\left(\left(\roman R_{a_t}\right)_* v_t'\right)
=
\roman{ad}_{a_t^{-1}}
\vartheta\left(v_t'\right).
}$
Further,
by construction,
$
\omega\left(v_t a_t'\right)
$
equals
${a^{-1}_t a_t' \in g = \roman T_eG}$.
Consequently the
fact that $\widehat u_t'$
is horizontal for the connection
$A+s \vartheta$ entails that
$a_t$ satisfies the differential equation
$$
0 =
a^{-1}_t a_t'
-
s\,
\roman{ad}_{a_t^{-1}}
\vartheta\left(v_t'\right) \in g
$$
in the Lie algebra $g$ of $G$ or, equivalently, the differential equation
$$
0 =
a_t'a^{-1}_t
-
s\,
\vartheta\left(v_t'\right) \in g.
$$
When we differentiate this equation with respect to $s$ we obtain
$$
0 =
\frac{\partial}{\partial s}\left(a_t'a^{-1}_t\right)
-
\vartheta\left(v_t'\right) \in g,
$$
that is,
$$
0 =
\left(\frac{\partial}{\partial s}a_t'\right) a^{-1}_t
+a_t' \left(\frac{\partial}{\partial s} a^{-1}_t\right)
-
\vartheta\left(v_t'\right) \in g.
\tag*
$$
\smallskip

Finally we observe that
by construction
the map $a$
is subject to the conditions
\linebreak
${
a(0,t) = e = a(s,0).
}$
In particular, for $s=0$,  do we have
$
a_t' = 0,
$
and hence, for $s=0$, the differential equation $(*)$
simplifies to
$$
0 =
\frac{\partial}{\partial s}a_t'
-
\vartheta\left(v_t'\right) \in g.
\tag**
$$
However, since
$a$ is defined on the product of two intervals,
we may interchange partial derivatives and obtain, for $s=0$,
the differential equation
$$
0 = \frac {\roman d}{\roman dt}
\frac{\partial a_t}{\partial s}
-
\vartheta\left(v_t'\right) \in g.
$$
{}From this we conclude
$$
\frac{\partial a}{\partial s}(0,t)
=
\int_0^t
\vartheta\left(v_{\tau}'\right) d\tau
=
\int_{u_{A,\widehat Q,t}}
\vartheta \ \in g. \qed
$$
\enddemo

\smallskip
By a {\it smooth\/} map $h$ on $\Cal A(\xi)$
with values in a smooth finite dimensional manifold
we mean henceforth
a continuous map $h$ whose restriction to an arbitrary smooth
finite dimensional submanifold of $\Cal A(\xi)$ is smooth in the ordinary
sense.
We can then still talk about the derivative of $h$:
for a point $A$ of $\Cal A(\xi)$, by the
{\it differential\/} or
{\it derivative\/}
$dh(A)$, {\it evaluated at a\/} 1-{\it form\/}
$\vartheta \in \Omega^1(M,\roman{ad}(\xi))$,
we mean the corresponding partial derivative.
\smallskip

For
a smooth closed path
$w \colon [0,b] \to M$, with starting point
$Q \in M$,
the {\it holonomy\/}
$\roman{Hol}_{w,\widehat Q}(A) \in G$
of $A$ {\it along\/} $w$
{\it with reference to\/}
$\widehat Q$
is defined by
\linebreak
${
u_{A,\widehat Q}(b) = \widehat Q\, \roman{Hol}_{w,\widehat Q}(A) \in P.
}$
For $y \in G$ we denote by
$L_y$ the operation of left translation
from $g$ to $\roman T_yG$.

\proclaim{Corollary 2.5}
When $u$ is closed
the
holonomy along $u$ furnishes a smooth
map
$\roman{Hol}_{u,\widehat Q}$
from
$\Cal A(\xi)$
to $G$.
Moreover, at a
connection $A$,
with $y={\roman{Hol}_{u,\widehat Q}(A)}\in G$,
the differential
$
d\roman{Hol}_{u,\widehat Q}(A)
\colon
\roman T_A\Cal A(\xi)
\to
\roman T_y G
$
assigns to a smooth 1-form
$
\vartheta \in
\roman T_A\Cal A(\xi)
=\Omega^1(M,\roman{ad}(\xi))
$
the value
${
L_y \int_{u_{A,\widehat Q}} \vartheta
\in
\roman T_y G}$.
Finally, this map
is invariant in the sense that, given a gauge transformation
$\gamma$, whatever smooth connection $A$,
${
\roman{Hol}_{u,\widehat Q}
(\gamma A)
=
x_{\gamma}
\roman{Hol}_{u,\widehat Q}
(A)
x_{\gamma}^{-1}}$.
(See Section 1 for the notation $x_{\gamma}$.)
\endproclaim

\demo{Proof}
Let $\vartheta \in \Omega^1(M,\roman{ad}(\xi))=\roman T_A\Cal A(\xi)$
and
$Y= \int_{u_{A,\widehat Q}} \vartheta \in g$.
By (2.1),
an explicit formula for the partial derivative
${
\frac{\partial U_b}{\partial \vartheta}(A)
=
dU_b(A)(\vartheta)}$
of the map
$U_b$ from $\Cal A(\xi)$ to  $P$ which assigns
$u_{A,\widehat Q}(b) \in P$ to $A \in \Cal A(\xi)$
is given by
$$
\frac{\partial U_b}{\partial \vartheta}(A)
=
\frac {\roman{d}}{\roman{dt}}
\left(u_{A,\widehat Q}(b)
 \roman{exp}t Y
\right)\big|_{t=0}\in
\roman T_{u_{A,\widehat Q}(b)}P.
$$
The derivative
$
\roman T_aG
@>>>
\roman T_{\widehat Q a}P
$
at $a \in G$
of the smooth map
from $G$ to $P$
which assigns to $a\in G$
the point
$\widehat Q a\in P$
may be described by the assignment
to
$L_a Z$
of
$\widehat Q L_a Z$,
for $Z \in g$. Hence
$$
\frac {\roman{d}}{\roman{dt}}
\left(u_{A,\widehat Q}(b)
 \roman{exp}t Y
\right)\big|_{t=0}
=
\frac {\roman{d}}{\roman{dt}}
\left(\widehat Q\, \roman{Hol}_{w,\widehat Q}(A)
 \roman{exp}t Y
\right)\big|_{t=0}
=
\widehat Q\,L_y \int_{u_{A,\widehat Q}} \vartheta. \qed
$$
\enddemo
\smallskip
Now we
pick
smooth closed curves
$w_1,\dots,w_n$
in $M$ starting at $Q$ and
representing  a set of generators
$x_1,\dots,x_n$
of the fundamental group $\pi= \pi_1(M,Q)$;
we write
$\bold w =(w_1,\dots,w_n)$ and
denote by $F$  the free group on
$x_1,\dots, x_n$.
The assignment to
a connection $A$ of
$\left(\roman{Hol}_{w_1,\widehat Q}(A),\dots,
\roman{Hol}_{w_n,\widehat Q}(A)\right) \in G^n$
yields a map
$$
\rho =
\roman{Hol}_{\bold w,\widehat Q}
\colon
\Cal A(\xi)
@>>>
G^n
\tag2.6
$$
which, in view of (2.1), is smooth
in the sense that its restriction to an arbitrary
smooth finite dimensional submanifold of $\Cal A(\xi)$ is smooth.
We refer to $\rho$ as
{\it Wilson loop mapping\/} since
its composite with
a smooth $G$-invariant function
on $\roman{Hom}(F,G)$
yields a smooth
$\Cal G(\xi)$-invariant function on
$\Cal A(\xi)$
generalizing what is
called a (classical) {\it Wilson loop observable\/}
in the physics literature.
Here is an immediate consequence of (2.5).

\proclaim{Theorem 2.7}
At a
connection $A$,
with
$$
\rho(A) =
\left(\roman{Hol}_{w_1,\widehat Q}(A),\dots,
\roman{Hol}_{w_n,\widehat Q}(A)\right)
=
(y_1,\dots,y_n) \in G^n,
$$
the differential
$
d\rho(A)
\colon
\roman T_A\Cal A(\xi)
@>>>
\roman T_{\rho(A)} G^n
=
\roman T_{y_1} G
\times
\dots
\times
\roman T_{y_n} G
$
of {\rm (2.6)} is given by the assignment to
$\vartheta \in \Omega^1(M,\roman{ad}(\xi))= \roman T_A\Cal A(\xi)$
of
$$
I_{\bold w,A,\widehat Q}(\vartheta) =
\left(
L_{y_1} \int_{\widehat w_1} \vartheta,
\dots,
L_{y_n} \int_{\widehat w_n} \vartheta
\right) \in
\roman T_{y_1} G
\times
\dots
\times
\roman T_{y_n} G,
$$
where, with an abuse of notation,
for $1 \leq j \leq n$,
$\widehat w_j$
denotes the horizontal lift of $w_j$
with reference to $A$ and
$\widehat Q$.
\endproclaim

\medskip\noindent{\bf 3. The first main result} \smallskip\noindent
In this Section we introduce our algebras of smooth
functions and spell out the first main result of the paper.
We return to the circumstances of the
Introduction.
Thus $\Sigma$
denotes
a closed surface,
$G$ a compact not necessarily connected
Lie group,
with Lie algebra $g$,
$\xi \colon P \to \Sigma$
a principal $G$-bundle,
having a connected total space $P$,
and $Q \in \Sigma$
a chosen  base point.
Consider the standard presentation
$$
\Cal P  = \big\langle x_1,y_1,\dots, x_\ell,y_\ell; r\big\rangle,
\quad
r = \left[x_1,y_1\right] \cdot
\dots
\cdot
\left[x_\ell,y_\ell\right] ,
\tag3.1
$$
of the fundamental group
$\pi=\pi_1(\Sigma,Q)$,
the number $\ell$ being
the genus of  $\Sigma$;
we denote by $F$ the free group on
$x_1,y_1,\dots, x_\ell,y_\ell$ and by $N$ the normal closure
of $r$ in $F$.
The quotient group $\Gamma = F\big/ [F,N]$
yields the {\it universal central extension\/}
$$
0
@>>>
\bold Z
@>>>
\Gamma
@>>>
\pi
@>>>
1
\tag3.2
$$
of
$\pi$; cf. Section 6 of \cite\atibottw\  and Section 2 of \cite\topology.
The topology of the bundle
$\xi$ determines an element
$X_{\xi}$
of the Lie algebra
$z$
of the centre of $Z$ of $G$
which is a topological characteristic class of $\xi$; see
\cite\atibottw\
for the
case of a connected structure group $G$
and Section 1 of our paper
{}~\cite\topology\
for the general case.
The evaluation map
which assigns
$\left(\phi(x_1),\phi(y_1),\dots, \phi(x_\ell),\phi(y_\ell)\right)
\in G^{2\ell}$ to
$\phi \in \roman{Hom}(F,G)$
identifies
$\roman{Hom}(F,G)$
with $G^{2\ell}$.
Let
$\roman{H}_{\xi}(\Gamma,G)$
be the
subspace of $\roman{Hom}(F,G)$
consisting of homomorphisms
$\chi \in \roman{Hom}(F,G)$
such that
$$
\left[\chi(x_1),\chi(y_1)\right] \cdot
\dots
\cdot
\left[\chi(x_\ell),\chi(y_\ell)\right]
=
\roman{exp}(X_{\xi}) \in Z.
\tag3.3
$$
This space is
manifestly compact and hence
has only finitely many
connected components;
furthermore
it is a finite union of real algebraic sets,
which, in turn,  also implies
that it
has only finitely many
connected components
since this is true of
any  real algebraic set, cf.
\cite\whitnfou.
\smallskip

The values  of the
restriction of the
Wilson loop mapping (2.6)
to the subspace
$\Cal N(\xi)$ of central Yang-Mills connections
lie in
$\roman{H}_{\xi}(\Gamma,G)$;
we denote by
$\roman{Hom}_{\xi}(\Gamma,G)$
its image
in $\roman{H}_{\xi}(\Gamma,G)$;
it is a space of homomorphisms from $\Gamma$ to $G$.
A more intrinsic description of the resulting
surjection from
$\Cal N(\xi)$  onto $\roman{Hom}_{\xi}(\Gamma,G)$
may be found  in (3.8) of our paper \cite\topology.
The connected components
of
$\roman{Hom}_{\xi}(\Gamma,G)$
are parametrized by the points of the corresponding
$\pi_0$-orbit in
$\roman{Hom}(\pi,\pi_0)$,
where $\pi_0$ refers to the group of connected components of $G$;
in particular, when $G$ is connected,
$\roman{Hom}_{\xi}(\Gamma,G)$
is connected.
Let
$I_{\xi}$ denote the ideal
in the algebra $C^{\infty}\left(\roman{Hom}(F,G)\right)$
of smooth functions
on $\roman{Hom}(F,G)$
that vanish on the subspace
$\roman{Hom}_{\xi}(\Gamma,G)$
of $\roman{Hom}(F,G)$, and {\it define\/}
an algebra $C^{\infty}\left(\roman{Hom}_{\xi}(\Gamma,G)\right)$
of continuous functions on
$\roman{Hom}_{\xi}(\Gamma,G)$ by
$$
C^{\infty}\left(\roman{Hom}_{\xi}(\Gamma,G)\right)=
C^{\infty}\left(\roman{Hom}(F,G)\right)/I_{\xi} .
\tag3.4
$$
This algebra is often called that of
{\it Whitney smooth functions\/} on
$\roman{Hom}_{\xi}(\Gamma,G)$,
cf.~\cite\whitnone.
We note that here and henceforth spaces
may arise which are not necessarily connected.
When we talk about an algebra of continuous functions
on such a space we always mean
an algebra of continuous functions
on a connected component.
We do not indicate this explicitly, to avoid
an orgy of notation.
\smallskip
Let
$\roman{Rep}_{\xi}(\Gamma,G)=\roman{Hom}_{\xi}(\Gamma,G)\big/ G$.
We {\it define\/}
an algebra
$C^{\infty}\left(\roman{Rep}_{\xi}(\Gamma,G)\right)$
of continuous functions
on
$\roman{Rep}_{\xi}(\Gamma,G)$
by
$$
C^{\infty}\left(\roman{Rep}_{\xi}(\Gamma,G)\right)
= \left(C^{\infty}(\roman{Hom}(F,G))\right)^G\big/I_{\xi}^G,
\tag3.5
$$
that is, we take
that of smooth
$G$-invariant functions $\left(C^{\infty}(\roman{Hom}(F,G))\right)^G$
on $\roman{Hom}(F,G)$
modulo its ideal
$I_{\xi}^G$
of
functions
that vanish on
$\roman{Hom}_{\xi}(\Gamma,G)$.
By construction this is an algebra of functions
on
$\roman{Rep}_{\xi}(\Gamma,G)$ in an obvious fashion.
Since $G$ is compact, the canonical map
from ${C^{\infty}\left(\roman{Rep}_{\xi}(\Gamma,G)\right)}$
to
$\left(C^{\infty}\left(\roman{Hom}_{\xi}(\Gamma,G)\right)\right)^G$
is a bijection whence
$C^{\infty}\left(\roman{Rep}_{\xi}(\Gamma,G)\right)$
may as well be described
as the algebra of
$G$-invariant
{\it Whitney smooth functions\/} on
$\roman{Hom}_{\xi}(\Gamma,G)$.
Since we shall not need this fact we refrain from giving
the details here.
\smallskip
In the same vein,
denote by  $C^{\infty}(\Cal A(\xi))$
the algebra of smooth functions on $\Cal A(\xi)$
in the sense explained in Section 2 above;
we then
{\it define\/}
the algebra $C^{\infty}\left(\Cal N(\xi)\right)$
on $\Cal N(\xi)$ as
the quotient algebra $C^{\infty}(\Cal A(\xi))/J_{\xi}$,
where
$J_{\xi}$ refers to the ideal of functions
in $C^{\infty}(\Cal A(\xi))$
that vanish on the subspace
$\Cal N(\xi)$ of $\Cal A(\xi)$, and
we {\it define\/}
an algebra $C^{\infty}\left(N(\xi)\right)$ of continuous functions
on the moduli space $N(\xi) =
\Cal N(\xi)\big/\Cal G(\xi)$
of central Yang-Mills connections
by
$$
C^{\infty}\left(N(\xi)\right)
= \left(C^{\infty}(\Cal A(\xi))\right)^{\Cal G(\xi)}\big/I_{\xi}^{\Cal G(\xi)},
\tag3.6
$$
that is, we take the algebra
of smooth ${\Cal G(\xi)}$-invariant functions
$\left(C^{\infty}(\Cal A(\xi))\right)^{\Cal G(\xi)}$
on $\Cal A(\xi)$
modulo its ideal
$I_{\xi}^{\Cal G(\xi)}$
of
functions
that vanish on
$\Cal N(\xi)$.
By construction, this is an algebra of functions
on
$N(\xi)$, in an obvious fashion.
\smallskip
The decomposition of $N(\xi)$
into connected components of
orbit types of classes of central Yang-Mills connections
endows $N(\xi)$
with a structure of a {\it decomposed space\/},
in fact,  see \cite\singulat~(1.2),
with that of a {\it stratified space\/}.
The {\it pieces\/} are
smooth manifolds,
parametrized by
conjugacy classes
$(K)$ of subgroups $K$
of $G$;
the piece
$N_{(K)}(\xi)$
corresponding to $(K)$
consists of classes
$[A]$ of central Yang-Mills connections
$A$
having
stabilizer
$Z_A \subseteq \Cal G(\xi)$
whose image in $G$ under (1.1) is conjugate to
$K$.
\smallskip

We now pick smooth closed paths
$u_1,v_1,\dots,u_\ell,v_\ell$
in $\Sigma$
representing the generators
$x_1,y_1,\dots,x_\ell,y_\ell$,
so that the standard cell decomposition of
$\Sigma$ with a single 2-cell $e$
corresponding to $r$ results,
and, furthermore,
a base point
$\widehat Q \in P$ so that
$\xi(\widehat Q) = Q \in \Sigma$.
Then the
Wilson loop mapping
$\rho$ from  $\Cal A(\xi)$ to $\roman{Hom}(F,G)$
with reference to these data,
cf. (2.6),
induces a homeomorphism
$$
\rho_{\flat} \colon N(\xi) @>>> \roman{Rep}_{\xi}(\Gamma,G);
\tag3.7
$$
it
coincides with the map
given in (3.8.2) of our paper \cite\topology.
By an abuse of language, we refer
to
$\rho_{\flat}$
as
{\it Wilson loop mapping\/} as well.
It is independent of the choices made to define $\rho$.
\smallskip

The decomposition
of $\roman{Rep}_{\xi}(\Gamma,G)$
into connected components of
orbit types of
representations
has as well pieces  parametrized by
conjugacy classes
$(K)$ of subgroups
of $G$;
the piece
$\roman R_{(K)}(\xi)$
corresponding to $(K)$
consists of classes
$[\phi]$ of
homomorphisms
$\phi$ from $\Gamma$  to $G$
having
stabilizer
$Z_{\phi} \subseteq G$
conjugate to
$K$.
The Wilson loop mapping
$\rho_{\flat}$
is manifestly compatible with
the decompositions
since (1.1)
identifies the stabilizer
$Z_A$
of a connection $A$ with
the stabilizer $Z_{\rho (A)}$
of $\rho (A) \in \roman{Hom}(F,G)$,
cf. e.~g.
\cite\singulat~(2.4).
Consequently the
Wilson loop mapping, restricted to a piece
$N_{(K)}(\xi)$ of $N(\xi)$,
is a homeomorphism onto
the corresponding piece
$\roman R_{(K)}(\xi)$
of the decomposition of $\roman{Rep}_{\xi}(\Gamma,G)$.
In particular,
each connected component of
a piece
$\roman R_{(K)}(\xi)$
of the decomposition
of $\roman{Rep}_{\xi}(\Gamma,G)$
into $G$-orbit types
inherits a  structure of a smooth manifold
from the corresponding stratum of
$N(\xi)$
in such a way that  this decomposition
of $\roman{Rep}_{\xi}(\Gamma,G)$ is as well a
stratification.
\smallskip
Given smooth spaces $(X,C^{\infty}(X))$
and $(Y,C^{\infty}(Y))$,
a map $\phi\colon X \to Y$
is said to be {\it smooth\/}
provided
for every
$f \in C^{\infty}(Y)$ the composite
$f \circ \phi$
is
a smooth function on $X$, that is, lies
in $C^{\infty}(X)$.
The usual notion of diffeomorphism carries over as well:
A smooth homeomorphism
is a {\it diffeomorphism\/}
provided its inverse map is also smooth.
Here is the {\it first main result\/} of the paper.

\proclaim {Theorem 3.8}
With reference to the
decompositions into connected components of orbit types,
the algebras
$C^{\infty}(N(\xi))$
and
$C^{\infty}(\roman{Rep}_{\xi}(\Gamma,G))$
yield smooth structures
on $N(\xi)$
and $\roman{Rep}_{\xi}(\Gamma,G)$, respectively,
and the Wilson loop mapping
$\rho_{\flat}$ from  $N(\xi)$
to $\roman{Rep}_{\xi}(\Gamma,G)$
is
a diffeomorphism
of smooth spaces.
\endproclaim

\demo{Remarks about the proof}
The restriction
of a
function
in $C^{\infty}\left(N(\xi)\right)$
to a piece
is a smooth function in the ordinary sense, and the same is true
of
$C^{\infty}(\roman{Rep}_{\xi}(\Gamma,G))$.
This is a consequence of the fact that the restriction of a smooth function
to a smooth submanifold is a smooth function on the submanifold.
A more formal proof
will be given in Section 6 below.
Hence the algebras $C^{\infty}(N(\xi))$
and
$C^{\infty}(\roman{Rep}_{\xi}(\Gamma,G))$
furnish
smooth structures as asserted.
Smoothness
of the map $\rho_{\flat}$
follows at once from the facts that
the Wilson loop mapping
$\rho$ from $\Cal A(\xi)$ to $\roman{Hom}(F,G)$
is smooth
and $\Cal G(\xi)$-invariant,
cf. (2.7),
where
$\Cal G(\xi)$
acts on
$\roman{Hom}(F,G)$
through the projection (1.1).
Moreover
$\rho^*$ is manifestly injective
since $\rho_{\flat}$ is a homeomorphism
and hence identifies the algebras of continuous functions
on these spaces.
The surjectivity of
$\rho^*$
will be established
in
Section 6 below by a partition of unity argument. \qed
\enddemo

\smallskip
Notice that a priori
the smooth structure $C^{\infty}(\roman{Rep}_{\xi}(\Gamma,G))$
depends on the choice
of presentation of $\pi$
but {\it not\/}
on the chosen Riemannian metric on $\Sigma$
while
the space
$N(\xi)$ and hence a fortiori its
smooth structure
$C^{\infty}(N(\xi))$
depend
on the chosen Riemannian metric on $\Sigma$
but {\it not\/}
on the choice
of presentation of $\pi$.
Theorem 3.8 implies that the smooth structure on
$\roman{Rep}_{\xi}(\Gamma,G)$ does {\it not\/}
depend on the choice of presentation.
Furthermore, a diffeomorphism
$\phi$ of $\Sigma$
preserving $\xi$
will induce a commutative diagram
$$
\CD
N(\xi) @>>>
\roman{Rep}_{\xi}(\Gamma,G)
\\
@V{\phi^{\sharp}}VV
@V{\phi^{\flat}}VV
\\
\widetilde N(\xi) @>>>
\roman{Rep}_{\xi}(\Gamma,G)
\endCD
$$
of diffeomorphisms of smooth spaces,
where
$\widetilde N(\xi)$
denotes the moduli space of central Yang-Mills connections
for the image under $\phi$ of the
chosen Riemannian metric on $\Sigma$.
We hope to return to this issue at another occasion.

\medskip\noindent{\bf 4. The twisted integration mapping in de Rham theory}
\smallskip\noindent
In the present Section
we work out a precise description of
the
twisted integration mapping
tailored to our purposes.
\smallskip

Consider a
principal $G$-bundle
$\xi \colon P \to M$
over an  {\it arbitrary\/}
smooth
connected finite dimensional
manifold $M$
having connected total space $P$.
As before we
pick a base point
$Q$ of $M$ and
a pre-image
$\widehat Q \in P$
of $Q$.
Given a flat connection $A$ on $\xi$,
the holonomy representation
$\phi =\rho(A)$
of $\pi = \pi_1(M,Q)$ in  $G$
induces  a  structure  of a $\pi$-module on
$g$ through
the adjoint action,
and it is folk lore that
the cohomology
$\roman H_A^*(M,\roman{ad}(\xi))$
is isomorphic to
the  cohomology
of $M$
with the appropriate local coefficients, cf. e.~g.
VII.7.3 on p. 107 of
\cite\raghuboo.
We need a more precise description
of a somewhat more general result, to be spelled out below.
\smallskip

Consider the universal covering
$\widetilde M \to M$
of $M$; we suppose that things have been set up in such a way
that
$\pi$
acts on the {\it right\/} of
$\widetilde M$, and we
pick a pre-image
$\widetilde Q \in \widetilde M$
of $Q$.

\proclaim{Proposition 4.1}
Every  smooth flat connection $A$ on $\xi$
determines a unique
smooth map
$\sigma=\sigma_{A,\widehat Q,\widetilde Q}$
from
$\widetilde M$  to $P$
which, with respect to the corresponding holonomy
representation
$\rho(A)$ of
$\pi$ in $G$,
furnishes a morphism
of (right) principal bundles over $M$.
\endproclaim

\demo{Proof}
This is established by an argument of the kind for the
{\it Reduction theorem\/} in II.7.1 of
\cite\kobanomi; for later reference we sketch the construction
of $\sigma$:
Given $T \in \widetilde M$,
let $\widetilde w$ be a smooth  path in
$\widetilde M$,
necessarily horizontal,
joining
$\widetilde Q$ and $T$,
let $w$ be the path in $M$ obtained by projecting
$\widetilde w$ into $M$,
and let
$\widehat w$
be the unique  lift
of $w$
that is horizontal
for $A$ and has starting point
$\widehat Q$;
then the value $\sigma (T)$ is defined as the end point
of
$\widehat w$.
Since $A$ is flat,
the value $\sigma (T)$
does not depend on the choice
of
$\widetilde w$. \qed
\enddemo

Let
$\zeta \colon E \to M$
be a smooth vector bundle
associated to
$\xi$
and the finite dimensional real representation
$V$ of $G$. Then
$\Omega^*(M,\zeta)$
amounts to the
$G$-invariant horizontal forms
in
$\Omega^*(P,V)$ and
the operator $d_A$ of covariant derivative of a {\it flat\/}
connection $A$ is a differential on
$\Omega^*(M,\zeta)$.
The following is immediate.

\proclaim{Corollary 4.2}
For every flat connection $A$,
the map
from
$\Omega^*(P,V)$  to $\Omega^*(\widetilde M,V)$
induced by $\sigma_{A,\widehat Q,\widetilde Q}$,
cf. {\rm (4.1)},
passes to an isomorphism
$
\sigma_{A,\widehat Q,\widetilde Q}^*
$
of chain complexes
from
$(\Omega^*(M,\zeta),d_A)$
onto
the subcomplex
$\left(\Omega^*(\widetilde M,V),d\right)^{\pi}$
of $\pi$-invariant $V$-valued forms
on $\widetilde M$,
the necessary $\pi$-module structure
on $V$ coming from the holonomy
$\pi \to G$ of $A$ combined with the $G$-action on $V$. \qed
\endproclaim

Given a homomorphism $\phi$ from $\pi$  to $G$
and a representation $V$ of $G$,
we write
$
\left(C^*(M,V_{\phi}),d\right)
$
for
the subcomplex
of $\pi$-invariant $V$-valued cellular cochains
on $\widetilde M$
and we denote by
$\roman H^*(M,V_{\phi})$
the resulting $\pi$-{\it equivariant\/}
cohomology of
$\widetilde M$ with values
in $V$.
It is naturally isomorphic to
the cohomology of $M$
with {\it local coefficients\/}
determined by
$\phi$ and the representation of
$G$ on $V$.
The usual integration mapping
${
\left(\Omega^*(\widetilde M,V),d\right)
@>>>
\left(C^*(\widetilde M,V),d\right)
}$
from the de Rham complex to that of
usual cellular
cochains is compatible with the $\pi$-actions.
Taking
invariants and combining it with
$\sigma_{A,\widehat Q,\widetilde Q}^*$,
for a given flat connection $A$,
we obtain the chain mapping
$$
\left(\Omega^*(M,\zeta),d_A\right)
@>>>
\left(\Omega^*(\widetilde M,V),d\right)^{\pi}
@>>>
\left(C^*(M,V_{\rho(A)}),d\right).
\tag4.3
$$
Henceforth we refer to it  as the {\it twisted integration mapping\/} in
de Rham theory;
it
induces an isomorphism
from
$\roman H^*_A(M,\zeta)$
onto $\roman H^*(M,V_{\rho(A)})$
a special case of which
is the folk lore isomorphism
mentioned earlier.
\smallskip

Under our circumstances,
twisted integration furnishes
such an isomorphism
even for a central connection which is not necessarily flat,
in the following way:
Recall \cite\singula\ that a
smooth connection $A$ on $\xi$
is said to be {\it central\/}
provided its curvature
$K_A$
is a 2-form on $M$ with values
in the Lie algebra $z$ of the centre $Z$ of $G$.
To apply what is said above
to a central connection,
write $Z_e$ for the
connected component of the identity of $Z$, let
$G^{\sharp}=G\big /Z_e$, $P^{\sharp}=P\big /Z_e$,
and consider
the induced principal
$G^{\sharp}$-bundle
$\xi^{\sharp} \colon P^{\sharp}@>>>M$;
since the adjoint representation
of $G$ on $g$ factors through
a representation of $G^{\sharp}$
the bundle
$\xi^{\sharp}$ is still a principal one for
$\roman{ad}(\xi)$.
Consequently a  central connection $A$
on $\xi$ induces a flat connection
$A^{\sharp}$ on
$\xi^{\sharp}$;
the operator $d_A$ of covariant derivative
is then a differential on
$\Omega^*(M,\roman{ad}(\xi))$, and
we can apply what is said above to
the vector bundle $\zeta =\roman{ad}(\xi)$
and corresponding principal bundle
$\xi^{\sharp}$.
Maintaining the notation established in Section 2, we suppose that
the smooth closed curves
$w_1,\dots,w_n$
are the 1-cells
of a cell decomposition of $M$ with the single zero cell $Q$,
and we thus in particular  identify
the fundamental group $\pi_1(M^1,Q)$
of the 1-skeleton $M^1$ of $M$
with the free group $F$.
Let, then, $A$ be a central connection
on $\xi$, and let
$\phi=\rho(A) \colon F \to G$.
With reference to the image
of $\widehat Q$ in $P^{\sharp}$,
the homomorphism $\phi$
manifestly passes to the standard holonomy homomorphism
from $\pi$ to $G^{\sharp}$
for the resulting flat connection
$A^{\sharp}$ on $\xi^{\sharp}$.
Abusing notation somewhat,
we write
$g_{\phi}$ for the
Lie algebra $g$ together with the
$\pi$-module structure induced by
$\phi$ and hence
by $A$;
the resulting twisted integration mapping,
with target the corresponding {\it cellular\/}
cochains,
then looks like
$$
\left(\Omega^*(M,\roman{ad}(\xi)),d_A\right)
@>>>
\left(C^*_{\roman{cell}}(M,g_{\phi}),d\right)
\tag4.4
$$
and induces, in particular,
an isomorphism
$
\roman{Int}_A
$
from
$\roman H_A^{*}(M,\roman{ad}(\xi))$
onto
$\roman H^{*}(M,g_{\phi})$.
When $M$ is aspherical,
the complex of cellular chains
$\bold C^{\roman{cell}}(\widetilde M)$
of the universal cover
$\widetilde M$
with its right $\pi$-module
structure
is a free resolution of the ground ring in the category of right
$\pi$-modules;
when $M$ is not aspherical,
a free resolution $\bold P$ is obtained by adding
to $\bold C^{\roman{cell}}(\widetilde M)$
more generators
in degrees $\geq 2$.
Consequently,
whatever right $\pi$-module $U$, the canonical map
from
$\roman H^*(\pi,U)$
to $\roman H^{*}(M,U)$
is an isomorphism in degree 1
and we shall take
it to be the identity,
the first cohomology of $\pi$ being computed from
$\bold P$.
Thus the isomorphism induced by the twisted integration mapping
furnishes, in degree 1,
an isomorphism
$
\roman{Int}_A
$
from
$
\roman H_A^{1}(M,\roman{ad}(\xi))$
onto
$\roman H^{1}(\pi,g_{\phi})$
while, for aspherical $M$,
in arbitrary degree,
it yields an isomorphism
$
\roman{Int}_A
$
from
$\roman H_A^{*}(M,\roman{ad}(\xi))$
onto
$\roman H^{*}(\pi,g_{\phi})$.

\define\de{\partial}

%\beginsection 5. Representation spaces

\bigskip\noindent{\bf 5. Representation spaces} \medskip\noindent
It
remains to rework
and extend
the classical relationship
between the
infinitesimal structure
of
representation spaces
and group cohomology,
cf. \cite\raghuboo, \cite\weilone, \cite\weiltwo.
Some care is necessary here since
central connections which are not necessarily flat
will come into play later.
\smallskip
Let
$$
\Cal P  = \big\langle x_1,\dots, x_n;
r_1,\dots,r_m\big\rangle,
\tag5.1
$$
be a presentation of
a finitely presented group
$\pi$, and
write $F$ for the free group on
$x_1,\dots, x_n$, so that $\pi = F/N$,
where $N$ refers to the normal closure of $r_1,\dots,r_m$.
Recall that,
given an element $w \in F$,
over any ground ring $R$,
the {\it right\/} {\smc Fox} derivative
$\frac {\de w}{\de x_j} \in RF$
with respect to the variable $x_j,\, 1 \leq j \leq n$,
is given by the equation
$$
1-w = \sum_{j=1}^n (1-x_j)\frac {\de w}{\de x_j}  \in IF.
$$
Here as usual
$
IK =
\roman{ker}\left(\varepsilon
\colon
RK
\longrightarrow
R\right)
$
refers to the {\it augmentation ideal\/}
of a group $K$.
The usual description
of a principal bundle with structure group acting on
the {\it right\/}
forces us to use here {\it right\/}
Fox derivatives
which are less common than
{\it left\/} Fox derivatives.
The Fox calculus, applies to
the presentation $\Cal P$,
yields
the sequence
$$
\widehat{\roman R(\Cal P)}
\colon
RF
@<{\de_1^F}<<
RF\left[x_1,\dots,x_n\right]
@<{\de_2^F}<<
RF\left[r_1,\dots,r_m\right]
\tag5.2
$$
involving the free
right $RF$-modules
having $r_1,\dots,r_m$ and
$x_1,\dots,x_n$ as bases, respectively;
further, the operators $\de_*^F$
are given by certain
explicit formulas; we reproduce them only for the case $m=1$,
which is our primary case of interest, and we write $r$ instead of $r_1$:
$$
\aligned
\de_2^F
&=
\left[
\frac {\de r}{\de x_1},\cdots,\frac {\de r}{\de x_n}
\right]^{\roman t}
\colon
RF\left[r\right]
@>{}>>
RF\left[x_1,\dots,x_n\right]
\\
\de_1^F
&=
\left[
1-x_1,
\cdots,
1-x_n
\right]
\colon
RF\left[x_1,\dots,x_n\right]
@>>>
RF ,
\endaligned
$$
where ${}^{\roman{t}}$ refers to the transpose of a vector.
Modulo $N$,
(5.2) yields the beginning
$\roman R(\Cal P)$
of
a free resolution of the ground ring $R$, viewed as a trivial $R\pi$-module,
in the category of right $R\pi$-modules;
the distinction between
$\roman R(\Cal P)$
and
$\widehat{\roman R(\Cal P)}$ will be  important in \cite\direct.
\smallskip

Given a right $RF$-module $U$,
with structure map
$\chi$ from $F$ to $\roman{Aut}(U)$,
application of the functor $\roman{Hom}_{RF}(-,U)$
to (5.2)
yields the sequence
$\roman{Hom}_{RF}(\widehat{\roman R(\Cal P)},U)$
which,
in view of the obvious identifications
$\roman{Hom}_{RF}(\widehat{\roman R_0(\Cal P)},U) =  U$,
$\roman{Hom}_{RF}(\widehat{\roman R_1(\Cal P)},U) =  U^n$,
$\roman{Hom}_{RF}(\widehat{\roman R_2(\Cal P)},U) =  U^m$,
looks like
$$
\roman{Hom}_{RF}(\widehat{\roman R(\Cal P)},U)
\colon
U
@>{\delta^0_{\chi}}>>
U^n
@>{\delta^1_{\chi}}>>
U^m.
\tag5.3
$$
Here the operators
$\delta_{\chi}$ depend on the
$RF$-module structure on $U$
while the modules
$U^m,\, U^n,\, U$
depend only on the presentation
whence the notation.
When
$\chi$
factors through a right $R\pi$-module structure on $U$,
(5.3) is a cochain complex
$(C^*(\Cal P,U),\delta^*_{\chi})$
computing low dimensional cohomology groups of
$\pi$ with coefficients in $U$.
Further, the subgroup of 1-cocycles
$Z^1(\Cal P,U)
= \roman{ker}(\delta^1_{\chi})$
then depends only on
$\pi,\,g$, and $\chi$,
and not on a choice of presentation
(5.1), and we shall
therefore write
$Z^1(\pi,U)$
instead of
$Z^1(\Cal P,U)$.
\smallskip

Henceforth we take
$R=\bold R$, the reals, and
$U=g$, the Lie algebra of $G$,
viewed as a {\it right\/}
$G$-module
in the usual way.
The assignment to
$\left(\chi(x_1),\dots, \chi(x_n)\right)$
of
\linebreak
$\chi \in \roman{Hom}(F,G)$ identifies $\roman{Hom}(F,G)$ with $G^{n}$,
and that of the $m$-tuple
$$
(r_1(\chi x_1,\dots,\chi x_n),\dots,r_m(\chi x_1,\dots,\chi x_n))
$$
to
$\chi \in \roman{Hom}(F,G)$
yields a
smooth map
$\Phi$
from
$\roman{Hom}(F,G)$ to $G^m$.
Moreover, for every
$\chi\in \roman{Hom}(F,G)$, we denote by
$\omega_{\chi}$
the
smooth map
from $G$
to $\roman{Hom}(F,G)$
which assigns
$x^{-1} \chi x \in \roman{Hom}(F,G)$
to $x \in G$.
For later reference we reproduce
the tangent behaviour of
these maps:
\smallskip
Let $\chi$
be a homomorphism
from $F$ to $G$; we write
$g_{\chi}$
for the Lie algebra $g$, viewed as a
right $F$-module
via $\chi$ and the adjoint  representation.
The
homomorphism
$\chi$
being viewed as the point
${
\bold y=
\left(y_1,\dots,y_n\right) =\left(\chi(x_1),\dots,\chi(x_n)\right)}$
of $G^n$,
its operation
of {\it left translation\/}
$\roman L_{\chi}$
from
$g^n$ to
$\roman T_{\chi}\roman{Hom}(F,G)$
amounts to
$
\roman L_{y_1}
\times
\dots
\times
\roman L_{y_n}
$
from
$g^n$
to
$\roman T_{y_1} G\times\dots\times\roman T_{y_n} G$.
Accordingly, we write
$
\roman L_{\Phi(\chi)}
$
for the corresponding operation of left translation
from
$g^m$
to
$\roman T_{\Phi(\chi)}G^m
=
\roman T_{r_1(\bold y)}G \times \dots \times\roman T_{r_m(\bold y)}G.
$
The following is well known, cf.
\cite\goldmone,
\cite\raghuboo,
\cite\weilone,~\cite\weiltwo.

\proclaim{Proposition 5.4}
The tangent maps
$\roman T_e\omega_{\chi}$
and $\roman T_{\chi}\Phi$
and the operations of left translation
make commutative the diagram
$$
\CD
\roman T_eG
@>\roman T_e\omega_{\chi}>>
\roman T_{\chi} \roman{Hom}(F,G)
@>{\roman T_{\chi} \Phi}>>
\roman T_{\Phi (\chi)}G^m
\\
@A{\roman{Id}}AA
@A{\roman L_{\chi}}AA
@A{\roman L_{\Phi(\chi)}}AA
\\
g
@>>{\delta^0_{\chi}}>
g^n
@>>{\delta^1_{\chi}}>
g^m
\endCD
$$
where
$\delta^0_{\chi}$ and
$\delta^1_{\chi}$ refer to the corresponding operators
in
{\rm (5.3)},
for $U = g_{\chi}$.
\endproclaim

\smallskip

For a homomorphism $\chi$ from  $F$ to $G$
having the property that each $\chi(r_j)$
lies in the centre of $G$,
the Lie algebra $g$
inherits a structure of a right $\pi$-module
which we still denote by
$g_{\chi}$.

\proclaim{Corollary 5.5}
At a homomorphism
$\chi$ from  $F$ to $G$
having the property that each $\chi(r_j)$
lies in the centre of $G$,
left translation
$
\roman L_{\chi}
$
from
$C^1(\Cal P,g_{\chi})=g^n$
to
$\roman T_{\chi}\roman{Hom}(F,G)$
identifies
the subspace
$Z^1(\pi,g_{\chi})$
of 1-cocycles
with
the kernel of the tangent map
$\roman T_{\chi}\Phi$
from
$\roman T_{\chi} \roman{Hom}(F,G)$
to
$\roman T_{\Phi (\chi)}G^m$
and, moreover,
the subspace
$B^1(\pi,g_{\chi})$
of 1-coboundaries
with
the tangent space
$\roman T_{\chi}(G \chi) \subseteq  \roman T_{\chi}\roman{Hom}(F,G)$
to the $G$-orbit
$G \chi$
of $\chi$ in $\roman{Hom}(F,G)$.
\endproclaim

\proclaim{Proposition 5.6}
For every  $\chi \in \roman{Hom}(F,G)$
having the property that each $\chi(r_j)$
lies in the centre of $G$,
for each $x \in G$, the vector space automorphism
$\roman{Ad}(x)$
of $g$
is an isomorphism
of right $\bold R\pi$-modules
from
$g_{\chi}$ to $g_{x\chi}$
and hence induces an isomorphism
$\roman{Ad}_{\flat}(x)$
from
$\roman H^1(\pi,g_{\chi})$ onto $\roman H^1(\pi,g_{x\chi})$.
\endproclaim

\demo{Proof} This is left to the reader. \qed
\enddemo

We now have the machinery in place to relate the
derivative of
the Wilson loop mapping (2.6) with twisted 1-cochains
and integration.
We suppose that
(5.1)
is the presentation $\Cal P$ of the fundamental group
$\pi =\pi_1(M,Q)$
having
generators
and relations
represented by the
smooth closed curves
$w_1,\dots,w_n$, cf. Sections 2 and 4 above,
and attaching maps of the 2-cells of the cell decomposition
of $M$, respectively.
\smallskip
Let  $A$
be a central connection
on
$\xi$,
let $\phi = \rho(A)\colon F \to G$
and, as before, write $g_{\phi}$
for the Lie algebra $g$, with the
$\pi$-module structure
induced by $\phi$.
Notice the
cellular 1-cochains
$C^1_{\roman{cell}}(M,g_{\phi})$
coincide with the 1-cochains
$C^1(\Cal P,g_{\phi})$
with reference to $\Cal P$, cf. (5.3).

\proclaim{Theorem 5.7}
The differential
$
d\rho(A)
\colon
\roman T_A\Cal A(\xi)
@>>>
\roman T_{\phi} \roman{Hom}(F,G)
$
of the Wilson loop mapping
$\rho$
from
$\Cal A(\xi)$ to
$\roman{Hom}(F,G)$
amounts to the composite
of the
twisted integration
mapping from
$\roman T_A\Cal A(\xi)=\Omega^1(M,\roman{ad}(\xi))$
to
$C^1(\Cal P,g_{\phi})$
with
left translation
$\roman L_{\phi}$
from
$C^1(\Cal P,g_{\phi}) = g^n$ to
$\roman T_{\phi} \roman{Hom}(F,G)$.
\endproclaim

\demo{Proof}
In view of
what was said about
the map from
$\widetilde M$ to $P$ in the proof of (4.1) and,
furthermore, in view of the description
(4.4) of the twisted integration mapping,
the statement
follows at once from (2.7)
and the fact
that
the cellular 1-cochains
$C^1_{\roman{cell}}(M,g_{\phi})$
coincide with the 1-cochains
$C^1(\Cal P,g_{\phi})$. \qed
\enddemo

\beginsection 6. Reduction of the smooth structures to  the local models

%\medskip\noindent{\bf 6. Reduction of the smooth structures to  the local
%models}\smallskip\noindent
We return to the situation of the Introduction.
For intelligibility we assemble at first a number of facts
established
in our papers ~\cite\singula~---~\cite\topology.
\smallskip

The orthogonal structure
on $g$
combined
with
the usual wedge product of forms
$\wedge$
and integration
induces a non-degenerate
bilinear pairing
$(\cdot,\cdot)$ between
$\Omega^*(\Sigma,\roman{ad}(\xi))$
and
$\Omega^{2-*}(\Sigma,\roman{ad}(\xi))$
given by
$
(\zeta, \lambda) =
\int_\Sigma \zeta \wedge\lambda .
$
In particular,
this furnishes a
{\it weakly\/} symplectic structure
$\sigma$ on
$\Omega^1(\Sigma,\roman{ad}(\xi))$
and hence one on $\Cal A(\xi)$, cf.
\cite\atibottw, \cite\singula~(1.1).
Furthermore,
the space $\roman \Omega^2(\Sigma,\roman{ad}(\xi))$
of 2-forms
being identified with
the dual of
$\roman \Omega^0(\Sigma,\roman{ad}(\xi))$  via
$(\cdot,\cdot)$,
the assignment
to a connection $A$ of its curvature $K_A$ yields
a momentum mapping
$J$ from
$\Cal A(\xi)$ to
$\Omega^2(\Sigma,\roman{ad}(\xi))$,
for the
action of the group $\Cal G(\xi)$ of gauge transformations
on $\Cal A(\xi)$, cf. \cite\atibottw.
\smallskip
Let $A$ be a central Yang-Mills connection, fixed until
further notice.
Its operator
${
d_A
\colon
\roman \Omega^*(\Sigma,\roman{ad}(\xi))
@>>>
\roman \Omega^{*+1}(\Sigma,\roman{ad}(\xi))
}$
of covariant derivative is a differential.
Hence the cohomology
$\roman H^*_A=\roman H^*_A(\Sigma,\roman{ad}(\xi))$ is defined.
The Lie bracket on $g$ induces a graded Lie algebra structure
$
[\cdot,\cdot]_A
$
on $\roman H^*_A$
and the orthogonal structure
on $g$
together with $(\cdot,\cdot)$
a non-degenerate graded bilinear  pairing
$(\cdot,\cdot)_A$
between
$\roman H^*_A$
and
$\roman H^{2-*}_A$.
In particular,
the latter
identifies
$\roman H^{2}_A$
with the dual of
the Lie algebra
$\roman H^{0}_A = z_A$
of the stabilizer
$Z_A \subseteq \Cal G(\xi)$
of $A$, and
the constituent
of $(\cdot,\cdot)_A$
in degree 1
is a symplectic structure
$\sigma_A$ on
$\roman H^{1}_A$.
Moreover
the assignment
to $\eta \in \roman H^1_A$ of
$\Theta_A(\eta)
=
\frac 12 [\eta,\eta]_A$
yields
a momentum
mapping
$\Theta_A$ from
$\roman H^1_A$ to $\roman H^2_A$
for the
$Z_A$-action on
$\roman H^1_A$,
cf. \cite\singula~(1.2.5),
in fact, the {\it unique\/} one
with $\Theta_A(0)=0$.
Write
$\roman H_A$
for its reduced space.
By~\cite\singula\ (2.32),
the reduced space $\roman H_A$
is a local model for
$N(\xi)$ near $[A]$ in the sense that
the data induce a
homeomorphism of a neighborhood
of $[0] \in \roman H_A$
onto a neighborhood of $[A]$
in $N(\xi)$. Our aim is to show
that
$\roman H_A$
is a local model near $[A]$
for
{\it all\/} the structure of interest to us.
To this end we observe first that
$\roman H_A$
inherits an obvious smooth structure
which we explain
under more general circumstances:
\smallskip
Let $M$ be a (finite dimensional)
symplectic
manifold, with a hamiltonian action of a compact Lie group
$K$ and momentum mapping
$\mu$ from $M$ to $k^*$, and let
$V = \mu^{-1}(0)$ denote its zero locus,
so that the reduced space looks like
$M_{\roman{red}} = V/K$.
With respect to the decomposition
into connected components of
orbit types,
the algebra
of Whitney smooth functions
$$
C^{\infty}(V)= C^{\infty}(M) \big/ I_V,
\tag6.1.1
$$
where $I_V$ refers to the ideal of functions that vanish on $V$,
endows
$V$
with a smooth structure;
likewise,
the algebra
$$
C^{\infty}(M_{\roman{red}})= C^{\infty}(M)^K \big/ (I_V\cap (C^{\infty}(M)^K))
\tag6.1.2
$$
yields
a smooth structure on the reduced space
in an obvious fashion,
where
$C^{\infty}(M)^K$
refers to the subalgebra
of $K$-invariant functions.
By construction,
$C^{\infty}(M_{\roman{red}})$
is an algebra of continuos functions on
$M_{\roman{red}}$.
In particular,
this construction,
applied to $M=\roman H^{1}_A$, $\mu =\Theta_A$, and
$K= Z_A$, yields the smooth space
$(\roman H_A,C^{\infty}(\roman H_A))$.
Our present aim is to show that
the latter is a local model
for $(N(\xi),C^{\infty}(N(\xi)))$
near $[A] \in N(\xi)$.
\smallskip
Let $(X,C^{\infty}(X))$
be a smooth space, and let
$Y$ be an open subset of $X$.
In order to avoid to have to talk about {\it sheaves\/}
of germs of smooth functions,
we {\it define\/}
a notion of {\it induced\/} smooth structure on $Y$ in the following way:
We shall say that a continuous function $f$ on $Y$ is {\it smooth\/}
if every point $y$ of $Y$ has an open neighborhood $U$
so that the restriction of $f$ to $U$ coincides with the restriction
to $U$
of a smooth function on $X$, that is, a member of
$C^{\infty}(X)$.
These smooth functions on $Y$ constitute an algebra
$C^{\infty}(Y)$
of continuous functions
on $Y$ which we refer to as its
{\it induced smooth structure\/}.
Notice the restriction map
from
$C^{\infty}(X)$
to
$C^{\infty}(Y)$
is {\it not\/} in general surjective.
When $X$ is a smooth manifold,
with its standard smooth structure, and $Y$
an open subset of $X$,
the algebra
$C^{\infty}(Y)$
is that
of smooth functions on $Y$ in the ordinary sense.

\proclaim{Theorem 6.2}
Near $[A] \in N(\xi)$,
the smooth space
$(\roman H_A,C^{\infty}(\roman H_A))$
is a local model for
$(N(\xi),C^{\infty}(N(\xi)))$.
More precisely, the choice of $A$
(in its class $[A]$)
induces a diffeomorphism
of an open neighborhood
$W_A$ of $[0] \in \roman H_A$
onto an open neighborhood
$U_A$
of
$[A] \in N(\xi)$,
where
$W_A$ and
$U_A$
are endowed with the induced
smooth structures
$C^{\infty}(W_A)$
and
$C^{\infty}(U_A)$,
respectively.
\endproclaim

\smallskip
To spell out
the representation space version of (6.2),
let
${
\phi  =\rho(A)
\colon
\Gamma
@>>>
G.
}$
Every $\psi \in \roman{Hom}_{\xi}(\Gamma,G)$
is manifestly of this form
and, given such a $\psi$,
a central Yang-Mills connection
on $\xi$
which is mapped to
$\psi$ under  $\rho$
is unique up to based gauge transformations;
see \cite\topology.
The same kind of structure
as that
denoted above by
$(\cdot,\cdot)_A$,
$\Theta_A$, and
$[\cdot,\cdot]_A$,
is available
on
$\roman H^{*}_\phi=\roman H^*(\pi,g_{\phi})$
and the twisted
integration mapping
from
$\roman H_A^*$
to
$\roman H^*_{\phi}$
identifies the respective structures.
In particular,
the Lie bracket on $g$ induces a graded Lie algebra structure
$[\cdot,\cdot]_{\phi}$
on $\roman H^*_{\phi}$.
Further,
the orthogonal structure on $g$ induces a
graded non-degenerate bilinear pairing
on $\roman H^*_{\phi}$
which in degree 1 amounts to a symplectic structure
$\sigma_{\phi}$ on $\roman H^1_{\phi}$,
and the assignment
to $\eta \in \roman H^1_{\phi}$
of
$\Theta_{\phi}(\eta)=\frac 12 [\eta,\eta]_{\phi}$
yields a momentum
mapping
$\Theta_{\phi}$ from
$\roman H^1_{\phi}$ to
$\roman H^2_{\phi}$,
for the
action
of the stabilizer
$Z_\phi \subseteq G$
of
$\phi \in \roman{Hom}_{\xi}(\Gamma,G)$
on
$\roman H^1_{\phi}$;
notice that the surjection (1.1)
passes to an isomorphism
from $Z_A$ to $Z_\phi$
identifying the stabilizers.
Moreover,
the construction (6.1.2),
applied to $M=\roman H^{1}_\phi$, $\mu =\Theta_\phi$, and
$K= Z_\phi$, yields the smooth space
$(\roman H_\phi,C^{\infty}(\roman H_\phi))$.

\proclaim{Theorem 6.3}
Near $[\phi] \in \roman{Rep}_{\xi}(\Gamma,G)$,
the smooth space
$(\roman H_\phi,C^{\infty}(\roman H_\phi))$
is a local model for the smooth space
$(\roman{Rep}_{\xi}(\Gamma,G),C^{\infty}(\roman{Rep}_{\xi}(\Gamma,G)))$.
More precisely, the choice of $\phi$
(in its class $[\phi]$)
induces a diffeomorphism
of an open neighborhood
$W_\phi$ of $[0] \in \roman H_\phi$
onto an open neighborhood
$U_\phi$
of
$[\phi] \in \roman{Rep}_{\xi}(\Gamma,G)$,
where
$W_\phi$ and
$U_\phi$
are endowed with the induced
smooth structures
$C^{\infty}(W_\phi)$
and
$C^{\infty}(U_\phi)$,
respectively.
\endproclaim

\proclaim{Addendum}
Under the circumstances of {\rm (6.2)} and {\rm (6.3)},
for suitable choices of the data,
twisted integration identifies the local models.
More precisely,
for a suitable choice of the data,
the twisted integration mapping
$\roman{Int_A}$
from
$\roman H^*_{A}$ to
$\roman H^*_{\phi}$
and the Wilson loop mapping
$\rho_{\flat}$
from
$N(\xi)$
to
$\roman{Rep}_{\xi}(\Gamma,G)$
fit into
a commutative diagram
$$
\CD
W_A
@>>>
U_A
\\
@V{\roman{Int_A}_{\sharp}}VV
@V{\rho_{\flat}|}VV
\\
W_\phi
@>>>
U_\phi
\endCD
$$
of diffeomorphisms between smooth spaces,
the four spaces being endowed with the smooth structures
mentioned earlier;
here $\roman{Int_A}_{\sharp}$
denotes the  map induced by twisted integration and
$\rho_{\flat}|$ the restriction of
the Wilson loop mapping to $U_A$, and the unlabelled horizontal arrows
are the maps coming into play in
{\rm (6.2)} and {\rm (6.3)}.
\endproclaim

\smallskip
The  proofs
of (6.2) and (6.3)
require some preparation.
Near $A$,
the pre-image
${
\Cal A_A
= J^{-1}\left(\Cal H_{A}^2(\Sigma,\roman{ad}(\xi))\right)
}$
of the space $\Cal H_{A}^2(\Sigma,\roman{ad}(\xi))$
of harmonic 2-forms
is a smooth $Z_A$-invariant submanifold of
$\Cal A(\xi)$, cf. \cite\singula,
and the operator $d_A$
gives rise to the exact sequence
$$
0
@>>>
\roman T_{A} \Cal A_A
@>>>
\roman T_{A} \Cal A(\xi)
@>{d_{A}}>>
\Omega^2(\Sigma,\roman{ad}(\xi))
@>>>
\roman H_{A}^2(\Sigma,\roman{ad}(\xi))
@>>>
0
\tag6.4
$$
of real vector spaces whence, in particular,
$\roman T_{A} \Cal A_A=Z_{A}^1(\Sigma,\roman{ad}(\xi))$,
the corresponding space of 1-cocycles;
here the tangent space $\roman T_{A} \Cal A(\xi)$
is identified with $\Omega^1(\Sigma,\roman{ad}(\xi))$ as usual.
Let $\Cal M_A$
be a  smooth
finite dimensional $Z_A$-invariant
submanifold
of
$\Cal A_A$
containing $A$,
of the kind
coming into play
in the proofs of
{}~\cite\singula\ (2.32) and
{}~\cite\singulat\ (1.2);
in particular,
${
\roman T_A \Cal M_A
=
\Cal H^1_A(\Sigma,\roman{ad}(\xi)),
}$
the subspace of harmonic 1-forms
in $\Omega^1(\Sigma,\roman{ad}(\xi))$;
in (6.11) below we shall pick
$\Cal M_A$
suitably.
We remind the reader that $\Cal N(\xi)\subseteq \Cal A(\xi)$
denotes the subspace of central Yang-Mills connections.
It is clear that the assignment
to a pair
$(\gamma, A)$
in
$\Cal G(\xi) \times \Cal A(\xi)$ of $\gamma (A)$
induces an injective
$\Cal G(\xi)$-invariant
immersion
$$
\Cal G(\xi) \times_{Z_A} \Cal M_A
@>>>  \Cal A(\xi)
\tag6.5
$$
identifying
$\Cal G(\xi) \times_{Z_A}\Cal M_A$
with a smooth
$\Cal G(\xi)$-invariant
codimension 0
submanifold of $\Cal A_A$
containing a
$\Cal G(\xi)$-invariant neighborhood of $A$
in $\Cal N(\xi)$.
In particular,
the derivative of this immersion at $A$
amounts to the inclusion
of $Z^1_A(\Sigma,\roman{ad}(\xi))$ into
$\Omega^1(\Sigma,\roman{ad}(\xi))$.
\smallskip
By ~\cite\singula~(2.18), the 2-form $\sigma$ on $\Cal A(\xi)$
passes to a symplectic structure
$\omega_A$ on
the smooth manifold $\Cal M_A$,
and  $J$
induces a momentum mapping
$\vartheta_A$
from
$\Cal M_A$
to
$\roman H^2_A(\Sigma,\roman{ad}(\xi))$
for the $Z_A$-action,
with $\vartheta(A) = 0$;
here $\roman H^2_A(\Sigma,\roman{ad}(\xi))$
is identified with the dual
$z_A^*$
of the Lie algebra $z_A= \roman H^0_A(\Sigma,\roman{ad}(\xi))$
as explained above;
see (2.21) in \cite\singula\ for details.
We now consider the
Marsden-Weinstein reduced space
$\Cal W_A = \vartheta_A^{-1}(0)\big/ Z_A$.
It is obvious that
(6.5)
induces an injection
$$
\Cal W_A
@>>>
N(\xi)
\tag6.6
$$
of $\Cal W_A$
into $N(\xi)$ identifying
$\Cal W_A$  with an open neighborhood
$U_A$
of
$[A]$ in $N(\xi)$,
and
in this way (6.6)
furnishes a {\it model\/}
of a neighborhood of
$[A]$ in $N(\xi)$.
Likewise,
the composite of
(6.6) with
the Wilson loop mapping
$\rho_{\flat}$ from $N(\xi)$ onto
$\roman {Rep}_{\xi}(\Gamma,G)$
is an injection
$$
\Cal W_A
@>>>
\roman {Rep}_{\xi}(\Gamma,G)
\tag6.7
$$
of $\Cal W_A$
into $\roman {Rep}_{\xi}(\Gamma,G)$ identifying
$\Cal W_A$  with an open neighborhood
$U_\phi$
of
$[\phi]$ in $\roman {Rep}_{\xi}(\Gamma,G)$
whence (6.7) furnishes a {\it model\/}
of a neighborhood of
$[\phi]$ in $\roman {Rep}_{\xi}(\Gamma,G)$.
With respect to the decompositions
into connected components of orbit types,
the embeddings (6.6) and (6.7)
are decomposition preserving.
The construction (6.1.2)
applied to $M=\Cal W_A$, $\mu =\vartheta_A$, and
$K= Z_A$, yields a smooth structure
$C^{\infty}(\Cal W_A)$
on $\Cal W_A$,
and
the embeddings (6.6) and (6.7) are smooth
since they preserve the decompositions into orbit types.
Let
$C^{\infty}(U_A)$
and
$C^{\infty}(U_\phi)$
be the
induced
smooth structures
on
$U_A$
and
$U_\phi$, respectively; it is obvious that
(6.6) and (6.7) induce smooth maps
$$
(\Cal W_A,C^{\infty}(\Cal W_A)
@>>>
(U_A,C^{\infty}(U_A))
\tag6.8
$$
and
$$
(\Cal W_A,C^{\infty}(\Cal W_A))
@>>>
(U_\phi,C^{\infty}(U_\phi)).
\tag6.9
$$
Moreover the Wilson loop mapping from
$N(\xi)$
to
$\roman {Rep}_{\xi}(\Gamma,G)$
passes to
a smooth map
$$
(U_A,C^{\infty}(U_A))
@>>>
(U_\phi,C^{\infty}(U_\phi))
\tag6.10
$$
in such a way that (6.9) is the composite of
(6.8) and (6.10).
Since each of (6.8), (6.9), (6.10)
are homeomorphisms between the underlying spaces,
the induced maps
$C^{\infty}(U_A)
\to
C^{\infty}(\Cal W_A)$
etc.
between the algebras
of smooth functions are injective.
We now show that they are surjective,
for a suitable choice of the data.
This will almost establish
the statements of (6.2) and (6.3),
except that
$\Cal W_A$
comes into play rather than
an open neighborhood
$W_A$ of $[0] \in \roman H_A$.
We proceed as follows:
\smallskip
The composite
$$
\Cal G(\xi) \times_{Z_A} \Cal M_A
@>>>
\roman {Hom}(F,G).
\tag6.11.1
$$
of
(6.5) with the Wilson loop mapping from
$\Cal A(\xi)$
to $\roman{Hom}(F,G)$
is
$\Cal G(\xi)$-invariant,
with respect to the
$\Cal G(\xi)$-action on
$\roman{Hom}(F,G)$
induced by (1.1) and, furthermore,
factors through the obvious  surjection
$$
\Cal G(\xi) \times_{Z_A} \Cal M_A
@>>>
G \times_{Z_A} \Cal M_A
\tag6.11.2
$$
and hence
passes to a smooth
$G$-invariant map
$$
G \times_{Z_A} \Cal M_A
@>>>
\roman {Hom}(F,G).
\tag6.11.3
$$

\proclaim{Proposition 6.11}
For a suitable choice of
$\Cal M_A$,
the map {\rm (6.11.3)}
is
a  smooth injective
$G$-invariant
immersion
identifying
$G\times_{Z_A}\Cal M_A$
with a
smooth
$G$-submanifold
of $\roman {Hom}(F,G)$
containing a
$G$-invariant neighborhood of $\phi$
in
$\roman {Hom}_{\xi}(\Gamma,G)$.
\endproclaim

To prepare for the proof,
we recall that
the tangent space
$\roman T_A\Cal M_A$
equals the space
$\Cal H^1_A(\Sigma,\roman{ad}(\xi))$
of harmonic 1-forms and
the tangent space
$\roman T_{(e,A)}
(\Cal G(\xi) \times_{Z_A} \Cal M_A)$
equals the space
$Z^1_A(\Sigma,\roman{ad}(\xi))$
of 1-cocycles; the latter, in turn,
decomposes into the direct sum of
$
B^1_A(\Sigma,\roman{ad}(\xi))
$
and
$\Cal H^1_A(\Sigma,\roman{ad}(\xi))$.
At
the point
$(e,A)$,
the tangent space
of
$G\times_{Z_A}\Cal M_A$
equals likewise
the direct sum of
$B^1(\pi,g_{\phi})$
and
$\Cal H^1_A(\Sigma,\roman{ad}(\xi))$,
and  the smooth map
(6.11.2)
has tangent map
$$
B^1_A(\Sigma,\roman{ad}(\xi))
\oplus
\Cal H^1_A(\Sigma,\roman{ad}(\xi))
@>{\left(\roman{Int}_A|,\roman{Id}\right)}>>
B^1(\pi,g_{\phi})
\oplus
\Cal H^1_A(\Sigma,\roman{ad}(\xi))
\tag6.11.4
$$
where $\roman{Int}_A|$ refers to the restriction
of the twisted integration mapping
$\roman {Int}_A$
from
$\Omega^*(\Sigma,\roman{ad}(\xi))$
to
$C^*(\Cal P,g_{\phi})$, cf. (4.4),
to the 1-coboundaries.
However,
the
restriction of the
twisted integration mapping
to the subspace
of 1-cocycles
$Z^1_A(\Sigma,\roman{ad}(\xi))$
amounts to
a surjection
of $Z^1_A(\Sigma,\roman{ad}(\xi))$
onto
$Z^1(\pi,g_{\phi})$, as
inspection of the commutative diagram
$$
\CD
0
@>>>
\roman H^0_A
@>>>
\Omega^0
@>>>
Z^1_A
@>>>
\roman H^1_A
@>>>
0
\\
@.
@VVV
@VVV
@VVV
@VVV
@.
\\
0
@>>>
\roman H^0_{\phi}
@>>>
C^0
@>>>
Z^1_{\phi}
@>>>
\roman H^1_{\phi}
@>>>
0
\endCD
$$
with the obvious unlabelled arrows
reveals,
where we have written
$\roman H^*_A=\roman H^*_A(\Sigma,\roman{ad}(\xi))$,
$\Omega^0=\Omega^0(\Sigma,\roman{ad}(\xi))$,
$Z^1_A=Z^1_A(\Sigma,\roman{ad}(\xi))$,
$\roman H^*_{\phi}=\roman H^*(\pi,g_{\phi})$,
$C^0=C^0(\Cal P,g_{\phi})$,
$Z^1_{\phi}=Z^1(\pi,g_{\phi})$
for short.
The diagram
has exact rows;
its outermost columns are isomorphisms;
and the arrow
from
$\Omega^0$
to
$C^0$
is manifestly surjective.
This implies that
(6.11.4)
is surjective.
In fact, write
$\Cal H^*(\pi,g_{\phi})$ for
the isomorphic image
in
$Z^*(\pi,g_{\phi})$
of the
subspace of harmonic forms
$\Cal H^*_A(\Sigma,\roman{ad}(\xi))$
in
$\Omega^*(\Sigma,\roman{ad}(\xi))$
under the twisted integration mapping
$\roman {Int}_A$
so that
the canonical epimorphism
from $Z^*(\pi,g_{\phi})$
onto $\roman H^*(\pi,g_{\phi})$
passes to an isomorphism
from
$
\Cal H^*(\pi,g_{\phi})
$
onto $\roman H^*(\pi,g_{\phi})$.
The direct sum
of $B^1(\pi,g_{\phi})$ and $\Cal H^1(\pi,g_{\phi})$
equals the space $Z^1(\pi,g_{\phi})$
of 1-cocycles,
and the
surjection
of $Z^1_A(\Sigma,\roman{ad}(\xi))$
onto
$Z^1(\pi,g_{\phi})$
factors through
the induced isomorphism
$(\roman{Id},\roman{Int}_A)$ from
$B^1(\pi,g_{\phi})\oplus\Cal H^1_A(\Sigma,\roman{ad}(\xi))$
onto
$B^1(\pi,g_{\phi})\oplus\Cal H^1(\pi,g_{\phi})$,
whence
(6.11.4)
is surjective.
Consequently
(6.11.2)
a submersion
near the point $(e,A)$.

\demo{Proof of {\rm (6.11)}}
The tangent map
of {\rm (6.11.3)}
at the point
$(e,A)$
is the composite of
\newline\noindent
{\rm (i)}
the isomorphism
$(\roman{Id},\roman{Int}_A)$ from
$B^1(\pi,g_{\phi})\oplus\Cal H^1_A(\Sigma,\roman{ad}(\xi))$
onto
$B^1(\pi,g_{\phi})\oplus\Cal H^1(\pi,g_{\phi})$,
\newline\noindent
{\rm (ii)}
the inclusion
of $B^1(\pi,g_{\phi})
\oplus
\Cal H^1(\pi,g_{\phi})
=Z^1(\pi,g_{\phi})$
into
$C^1(\Cal P,g_{\phi})$
and, finally,
\newline\noindent
{\rm (iii)}
left translation
$\roman L_{\phi}$ from
$C^1(\Cal P,g_{\phi})$ to $\roman T_{\phi}\roman{Hom}(F,G)$.
\newline\noindent
In fact,
in view (5.7),
the derivative of
(6.11.1)
at $A$
amounts to the twisted integration mapping
$\roman {Int}_A$
from
$\Omega^1(\Sigma,\roman{ad}(\xi))$
to
$C^1(\Cal P,g_{\phi})$, restricted to
the tangent space
$\roman T_A \Cal A_A=
Z_A^1(\Sigma,\roman{ad}(\xi)) \subseteq \Omega^1(\Sigma,\roman{ad}(\xi))$,
combined with
left translation
$\roman L_{\phi}$
from $C^1(\Cal P,g_{\phi})$ to
$\roman T_{\phi}\roman{Hom}(F,G)$.
However
it is manifest that
this tangent map
factors through
the map
from
$Z_A^1(\Sigma,\roman{ad}(\xi))$
onto
$Z^1(\pi,g_{\phi})
=
B^1(\pi,g_{\phi})
\oplus
\Cal H^1(\pi,g_{\phi})
$
induced by
the twisted integration mapping
and hence through (6.11.4).
Hence the tangent map
of {\rm (6.11.3)}
at the point
$(e,A)$
decomposes into the three pieces
(i) -- (iii) and is therefore
injective
since so is
the inclusion
of $Z^1(\pi,g_{\phi})$
into $C^1(\Cal P,g_{\phi})$.
This implies that
the smooth map {\rm (6.11.3)}
is an immersion near $(e,A)$;
hence, for a suitable choice of
$\Cal M_A$, it is injective.
\smallskip
Finally, since
$\Cal G(\xi) \times_{Z_A}\Cal M_A$,
viewed as
a smooth
$\Cal G(\xi)$-invariant
codimension 0
submanifold of $\Cal A_A$
via (6.5),
contains a
$\Cal G(\xi)$-invariant neighborhood of $A$
in $\Cal N(\xi)$,
and, furthermore, since
(6.11.2) is a submersion,
the image of
$G\times_{Z_A}\Cal M_A$
under (6.11.3)
contains a
$G$-invariant neighborhood of $\phi$
in
$\roman {Hom}_{\xi}(\Gamma,G)$ as asserted. \qed
\enddemo

Henceforth we assume that
the smooth manifold $\Cal M_A$
has been chosen in such a way that
(6.11.3) is injective.
This enables us to relate
the smooth structures of
$N(\xi)$ near $[A]$
and of $\roman {Rep}_{\xi}(\Gamma,G)$ near $[\phi]$
with that of
$\Cal W_A$
near $A$ by means of (6.11.3).
\smallskip
To verify surjectivity of
the induced
map from
$C^{\infty}(U_\phi)$
to
$C^{\infty}(\Cal W_A)$,
let $h\colon  \Cal W_A \to \bold R$
be a function in  $C^{\infty}(\Cal W_A)$.
Then there is a unique continuous function
$f$ on $U_\phi$ whose composite with
(6.9) equals $h$. We must show that
$f$ lies in $C^{\infty}(U_\phi)$.
In order to see this,
let $H$
be a smooth $Z_A$-invariant
function
on $\Cal M_A$
representing $h$.
Abusing notation, we denote its canonical extension
to a $G$-invariant function on
$G \times _{Z_A}\Cal M_A$
by $H$ as well.
The space $G \times _{Z_A}\Cal M_A$
being identified with a smooth
$G$-invariant submanifold of
$\roman{Hom}(F,G)$
via (6.11.3),
we must show that
$H$ extends locally to a $G$-invariant function
on
$\roman{Hom}(F,G)$.
However,
given
a homomorphism
$\psi$
from $\Gamma$ to $G$
in the image of
(6.11.3),
there is an open $G$-invariant neighborhood
$U$ of
$\psi$
in the image of
(6.11.3)
and a smooth $G$-invariant function
$\widetilde H$ on $\roman{Hom}(F,G)$
whose restriction to $U$ coincides with the restriction
of $H$ to $U$.
By construction,
$\widetilde H$
represents
a function
in $C^{\infty}(\roman {Rep}_{\xi}(\Gamma,G))$
and hence one
in
$C^{\infty}(U_\phi)$
which,
on a neighborhood of
$[\psi]$
in
$U_\phi$,
coincides with $f$.
Since $\psi$ is arbitrary,
this shows that
$f$
is smooth as asserted,
that is, lies in
$C^{\infty}(U_\phi)$.
Consequently (6.9), and hence (6.8) and (6.10),
are diffeomorphisms of smooth spaces.
\smallskip
To complete the proofs of (6.2) and (6.3)
we recall that,
by \cite\singula\ (2.31),
a suitable Kuranishi map
furnishes a
$Z_A$-equivariant
symplectomorphism
$
\Phi_A$
from
$\Cal M_A$
onto
a $Z_A$-invariant
ball
$B_A$
in $\roman H^1_A(\Sigma,\roman{ad}(\xi))$
about the origin,
cf. \cite\singula\ (2.29),
and this map
preserves the momentum mappings
$\Theta_A$ and
$\vartheta_A$.
Marsden-Weinstein reduction applied to
$B_A$
and $\Theta_A$,
restricted to $B_A$,
then yields the
open subspace $W_A$
of $\roman H_A$ we are looking for,
and the Kuranishi map
induces a homeomorphism
of a neighborhood of $[A]$ in $N(\xi)$
onto
$W_A$.
See~\cite\singula\ (2.32)
for details.
Moreover the construction (6.1.2)
applied to $M= B_A$, $\mu =\Theta_A$,
restricted to $B_A$,
and
$K= Z_A$, yields a smooth structure
$C^{\infty}(W_A)$
in such a way that
$\Phi_A$ induces
a diffeomorphism
from
$(\Cal W_A,C^{\infty}(\Cal W_A))$ onto
$(W_A,C^{\infty}(W_A))$.
Hence the
data
induce a diffeomorphism
of
$(W_A,C^{\infty}(W_A))$
onto
$(U_A,C^{\infty}(U_A))$.
This completes the proof of Theorems 6.2.
The same construction applies to the image
$B_\phi$ of
$B_A$
in
$\roman H^1_{\phi}$
under the twisted integration mapping $\roman{Int}_A$
from $\roman H^1_{A}$ to $\roman H^1_{\phi}$,
the momentum mapping
$\Theta_{\phi}$, and the stabilizer
$Z_\phi$ of
$\phi$; it yields
the open subspace $W_\phi$
of $\roman H_\phi$ we are looking for
and
a smooth structure
$C^{\infty}(W_\phi)$,
together with a diffeomorphism
of
$(W_\phi,C^{\infty}(W_\phi))$
onto
$(U_\phi,C^{\infty}(U_\phi))$.
This completes the proof of Theorems 6.3.
Moreover, the constructions have been carried out in such a way that
the statement of the Addendum is immediate. \qed
\smallskip
We now proceed towards the proof of
(3.8).
Henceforth $A$ will denote a central Yang-Mills connection
which is no longer fixed.
At first, we must show that
the restriction of a smooth function
on $N(\xi)$ and likewise on
$\roman {Rep}_{\xi}(\Gamma,G)$
to a stratum is a smooth function
on the stratum in the ordinary sense.
In view of (6.2) and (6.3),
it suffices to prove that,
under the circumstances of the construction
(6.1.2),
the restriction of a smooth
function in $C^{\infty}(M_{\roman{red}})$
to a piece
of
$M_{\roman{red}}$
is smooth
in the ordinary sense.
However this amounts to the fact that the
restriction
to a smooth submanifold of a smooth function
defined on a smooth manifold
is smooth on the submanifold.
\smallskip
As immediate consequence
of the Addendum to (6.3)
we see that
the Wilson loop mapping
from
$N(\xi)$ to
$\roman{Rep}_{\xi}(\Gamma,G)$
is {\it locally\/} a diffeomorphism.
To see that this is
{\it globally\/} so,
we
establish the existence of suitable partitions of unity.
We begin with the following
the proof of which is routine and therefore
left to the reader.

\proclaim{Lemma 6.12}
Let $W$ be a finite dimensional complex
representation of a compact Lie group $K$,
and let $B$ be an open $K$-invariant neighborhood of the origin.
Then there are open $K$-invariant  neighborhoods
$Q$ and $R$ of the origin with
$\overline Q \subseteq R$ and
$\overline R \subseteq B$, together with
a smooth $K$-invariant real-valued function
$H$ on $B$ with
$$
H\big | \overline Q = 1,
\quad
H\big | B \setminus R= 0.
$$
\endproclaim

Under the circumstances of (6.12), suppose the $K$-representation
is unitary, let $\mu$ denote its unique momentum mapping
from $W$ to $k^*$ having the value zero at the origin,
let $W_{\roman{red}}$ be its reduced space,
and let
$C^{\infty}(W_{\roman{red}})$ be the corresponding smooth structure (6.1.2).
Here is an immediate consequence of (6.12).

\proclaim{Corollary 6.13}
Let $P$ be an open neighborhood
in $W_{\roman{red}}$
of the class $[0]$ of the origin,
with its induced
smooth structure
$C^{\infty}(P)$.
Then there are open  neighborhoods
$Q$ and $R$
in $W_{\roman{red}}$
of  $[0]$ with
$\overline Q \subseteq R$ and
$\overline R \subseteq P$, together with
a smooth function $h \in C^{\infty}(P)$,
with
$$
h\big | \overline Q = 1,
\quad
h\big | P \setminus R= 0.
$$
\endproclaim

\proclaim{Corollary 6.14}
Given an arbitrary open neighborhood $U_{[A]}$ of
the point $[A]$ of $N(\xi)$,
there are open neighborhoods $Q_{[A]}$ and $R_{[A]}$ of $[A]$ in $N(\xi)$,
with $\overline Q_{[A]} \subseteq R_{[A]}$ and
$\overline R_{[A]} \subseteq U_{[A]}$,
together with a smooth function
$h_{[A]}$ on $N(\xi)$
with
$$
h_{[A]}\big | \overline Q_{[A]} = 1,
\quad
h_{[A]}\big | N(\xi) \setminus R_{[A]}= 0.
$$
Likewise, given an arbitrary open neighborhood $U_{[\phi]}$ of
the point $[\phi]$ of $\roman{Rep}_{\xi}(\Gamma,G)$,
there are open neighborhoods $Q_{[\phi]}$ and $R_{[\phi]}$ of $[\phi]$ in
$\roman{Rep}_{\xi}(\Gamma,G)$,
with $\overline Q_{[\phi]} \subseteq R_{[\phi]}$ and
$\overline R_{[\phi]} \subseteq U_{[\phi]}$,
together with a smooth function
$h_{[\phi]}$ on $\roman{Rep}_{\xi}(\Gamma,G)$
with
$$
h_{[\phi]}\big | \overline Q_{[\phi]} = 1,
\quad
h_{[\phi]}\big | \roman{Rep}_{\xi}(\Gamma,G) \setminus R_{[\phi]}= 0.
$$
When $[\phi] = \rho_{\flat} [A]$
and
$U_{[\phi]} = \rho_{\flat}(U_{[A]})$,
under the Wilson loop mapping
$\rho_{\flat}$ from $N(\xi)$ to
$\roman{Rep}_{\xi}(\Gamma,G)$,
things may be arranged in such a way that $\rho_{\flat}$
identifies
$Q_{[A]}$, $R_{[A]}$, and $h_{[A]}$
with respectively
$Q_{[\phi]}$, $R_{[\phi]}$, and $h_{[\phi]}$.
\endproclaim

\demo{Proof}
This is a consequence of (6.2), (6.3), its Addendum, and (6.13). \qed
\enddemo

For each
point $[A]$ of $N(\xi)$, pick
an injection of $W_A$ into $N(\xi)$
of the kind  coming into play in (6.2) above,
and write
$U_A \subseteq N(\xi)$
for the image of
$W_A$
in $N(\xi)$,
so that
$U_A$ is an open neighborhood of
$[A]$ in $N(\xi)$, as in (6.2);
we then write
$\phi = \rho(A)$ and
$U_\phi \subseteq
\roman{Rep}_{\xi}(\Gamma,G)$
for the image of
$U_A$
under the Wilson loop mapping, as in (6.3).
Here is our {\it third main result\/}.

\proclaim{Theorem 6.15}
There is a finite open cover of
$N(\xi)$ by open sets of the kind
$U_A$
together with a smooth partition of unity
subordinate to this cover.
Moreover,
there is a finite open cover of
$\roman{Rep}_{\xi}(\Gamma,G)$ by open sets of the kind
$U_\phi$
together with a smooth partition of unity
subordinate to this cover
in such a way that the Wilson loop mapping identifies
the covers and partitions of unity.
\endproclaim

\demo{Proof}
By (6.14),
for every central Yang-Mills connection
$A$,
there are open neighborhoods $Q_{[A]}$ and $R_{[A]}$ of $[A]$ in $N(\xi)$,
with $\overline Q_{[A]} \subseteq R_{[A]}$ and
$\overline R_{[A]} \subseteq U_{[A]}(\xi)$,
together with a smooth function
$h_{[A]}$ on $N(\xi)$
with
$$
h_{[A]}\big | \overline Q_{[A]} = 1,
\quad
h_{[A]}\big | N(\xi) \setminus R_{[A]}= 0.
$$
The subsets
$Q_{[A]}$ constitute an open cover of
$N(\xi)$.
Since $N(\xi)$ is compact, there is a finite subcover
$\{Q_1,\dots,Q_m\}$.
Each
$Q_\lambda$
lies in some
$U_\lambda$;
the corresponding family
$\{U_{\lambda}\}$
is the open cover of $N(\xi)$ we are aiming at.
Moreover, for each $\lambda$,
there is a function
$h_{\lambda} \in C^{\infty}(N(\xi))$
so that
$h_{\lambda}$ has the constant value 1 on
$\overline Q_{\lambda}$
and is zero outside an open neighborhood of
$\overline Q_{\lambda}$
in
$U_{\lambda}$.
Let $h = \sum h_{\lambda}$;
then $h \in C^{\infty}(N(\xi))$
and
$h[A] \geq 1$, whatever $[A] \in N(\xi)$.
The family
$\{e_{\lambda}\}$, where
$e_{\lambda} = \frac{h_{\lambda}}{h}$,
then furnishes the desired partition of unity.
\smallskip
The same kind of construction yields the asserted open cover and smooth
partition of unity for
$\roman{Rep}_{\xi}(\Gamma,G)$,
and the Wilson loop mapping identifies
the covers and partitions of unity. \qed
\enddemo

We can now complete
the proof of Theorem 3.8:
Let
$\{h_1,\dots,h_m\}$ be the partition of unity
subordinate to  the open cover
$\{U_1,\dots,U_m\}$
in (6.15).
Given $f \in
C^{\infty}\left(N(\xi)\right)$,
let
$f_{\lambda} = f h_{\lambda}$;
this is a smooth function, that is,
$f_{\lambda}\in C^{\infty}\left(N(\xi)\right)$.
By construction, each
$f_{\lambda}$
has a pre-image in
$C^{\infty}\left(\roman{Rep}_{\xi}(\Gamma,G)\right)$.
Consequently
$f$ has a pre-image in
$C^{\infty}\left(\roman{Rep}_{\xi}(\Gamma,G)\right)$
whence the map from
$C^{\infty}\left(\roman{Rep}_{\xi}(\Gamma,G)\right)$
to
$C^{\infty}\left(N(\xi)\right)$
induced by the Wilson loop mapping
is surjective.
This completes the proof of Theorem 3.8.

\medskip\noindent{\bf 7. Cohomology, Zariski tangent spaces,
and local semi-algebraicity}
\smallskip\noindent
In this Section we study the infinitesimal
structure of our spaces of interest.
\smallskip
Given a
smooth space $(X,C^{\infty}(X))$,
for each point $x \in X$,
the {\it ideal\/}
$\bold m_x$
of $x$ consists of all functions
in $C^{\infty}(X)$
vanishing at $x$;
as usual,
the space of {\it differentials\/}
$\Omega_x(X)$  at $x$ is the vector space
${
\Omega_x(X)
=\bold m_x
\big/\bold m^2_x,
}$
and
the {\it Zariski tangent space\/} $\roman T_x X$
is the dual space
${
\roman T_x X
=
\Omega_x(X)^*
=
\left(\bold m_x
\big/\bold m^2_x\right)^* .
}$
When $X$ is a smooth manifold near a point $x$
in the usual sense,
with standard smooth structure near $x$,
the  Zariski
tangent space boils down
to the usual smooth tangent space $\roman T_xX$
whence there is no risk of confusion in notation.
Here is another well known description
of the
Zariski tangent space:
Let $x \in X$ and
view $\bold R$ as a
$C^{\infty}(X)$-module,
written
$\bold R_x$,
by means of
the evaluation mapping
from $C^{\infty}(X)$ to $\bold R$
which assigns to a function
$f$ its value $f(x)$ at $x \in X$;
now
a {\it derivation\/}
at $x \in X$
is a linear map
$d$ from $C^{\infty}(X)$ to $\bold R$
satisfying the usual {\it Leibniz\/} rule
$$
d(fh) = (df)h(x) + f(x)dh .
$$
We denote the real vector space
of all
derivations
of
$C^{\infty}(X)$ in $\bold R_x$
by $\roman{Der}(C^{\infty}(X),\bold R_x)$.
For $x \in X$, the assignment
to $\phi \in \roman T_x X$
of the derivation
$d_{\phi}$ at $x$
given by
$d_{\phi} (f) = \phi( f - f_x)$
identifies
$\roman T_x X$
with
$\roman{Der}(C^{\infty}(X),\bold R_x)$;
here $f \in C^{\infty}(X)$ and
$f_x$ denotes the  function
having constant value
$f(x)$.
\smallskip
Given smooth spaces $(X,C^{\infty}(X))$,
$(Y,C^{\infty}(Y))$,
and a smooth  map $\phi$ from $X$  to $Y$,
the {\it derivative\/}
at a point $x \in X$
is the dual
${
d\phi_x
\colon
\roman T_x X
@>>>
\roman T_{\phi x} Y
}$
of the linear map
from
$\bold m_{\phi(x)}\big/\bold m^2_{\phi(x)}$
to
$\bold m_x\big/\bold m^2_x$
induced by $\phi$.
\smallskip
Let $(X,C^{\infty}(X))$
be a smooth space, and let
$U$ be an open subset of $X$.
We shall say that
a smooth function
$h$ on $X$
is a {\it bump \/} function
{\it with support in\/} $U$
if
there are open subsets $Q$ and $R$ of $X$
with $\overline Q \subseteq R$ and
$\overline R \subseteq U$,
so that
$$
h\big | \overline Q = 1,
\quad
h\big | X \setminus R= 0.
$$
Given
a point $x$ of $X$, we shall say that
$X$ {\it has smooth bump functions arbitrarily close to\/}
$x$ if for every open neighborhood $U$ of $x$ in $X$
there is a smooth bump function $h$ having the value
1 near $x$, with support in $U$.
{}From Corollary 6.14 above
we deduce at once the following.

\proclaim{Proposition 7.1} The spaces $N(\xi)$ and
$\roman {Rep}_{\xi}(\Gamma,G)$
have smooth bump functions arbitrarily close
to every point.
\endproclaim

Let $(X,C^{\infty}(X))$
be a smooth space
having smooth bump functions
arbitrarily close
to every point. We recall the following
well known
fact and reproduce a proof for completeness:

\proclaim{Proposition 7.2}
For every connected open subset $Y$, with
induced smooth structure
$C^{\infty}(Y)$,
the inclusion $j$ from $Y$ to $X$
induces an isomorphism
of Zariski tangent spaces for every $x \in Y$.
\endproclaim

\demo{Proof}
If $f$ is a smooth function which is constant on a neighborhood
$U$ of $x \in X$,
then $df$ is zero for every derivation
$d$ from $C^{\infty}(X)$ to $\bold R_x$.
In fact,
the differential of a constant function
(on $X$)
is zero, and hence
we may assume that $f$ has the value zero on $U$.
Given  a bump function $h$
with support in $U$, we then have
$$
0 = d(fh) = df h(x) + f(x) dh = df
$$
since $h(x) = 1$ and $f(x) = 0$.
\smallskip
In particular,
for every derivation
$d$ from $C^{\infty}(X)$ to $\bold R_x$,
the value $dh$ is zero for every bump function $h$
near $x \in X$.
Hence, given an arbitrary function
$f  \in C^{\infty}(X)$ and a bump function $h$ near $x$,
for every derivation
$d$ from $C^{\infty}(X)$ to $\bold R_x$,
we have
$$
d(fh) = (df) h(x) + f(x) dh = df.
$$
\smallskip
Let $x \in Y$, and let
$h$ be bump function on $X$ with
$h(y) =1$ near $x$ having support in $Y$.
Given a derivation
$d$ from
$C^{\infty}(X)$ to $\bold R_x$
and
$f \in C^{\infty}(Y)$,
the function $fh$ is defined on $X$,
and
$df=d(fh)$
extends $d$ to a
derivation from
$C^{\infty}(Y)$ to $\bold R_x$.
This shows the induced map from
$\roman{Der}(C^{\infty}(Y),\bold R_x)$ to
$\roman{Der}(C^{\infty}(X),\bold R_x)$
is surjective.
Moreover if
a derivation $d$
from
$C^{\infty}(Y)$ to $\bold R_x$
goes to zero in
$\roman{Der}(C^{\infty}(X),\bold R_x)$,
it must itself be zero since
$df = d(fh)$ for every $f$ and every bump function $h$. \qed
\enddemo

In view of (7.2),
there is no need for us to talk about {\it sheaves\/}
of germs of smooth functions
in order to define Zariski tangent spaces etc.
In fact, in view of (6.2), (6.3), and (7.1),
(7.2) entails at once the following:

\proclaim{Theorem 7.3} For every central Yang-Mills connection $A$,
the inclusion
of an open subspace
of the kind
$U_A$
into $N(\xi)$
induces an isomorphism
of Zariski tangent spaces
from
$\roman T_{[A]} U_A$
onto
$\roman T_{[A]} N(\xi)$.
Likewise, for every $\phi \in \roman{Hom}_{\xi}(F,G)$,
the inclusion
of an open subspace
of the kind
$U_\phi$
into
$\roman {Rep}_{\xi}(\Gamma,G)$
induces an isomorphism
of Zariski tangent spaces
from
$\roman T_{[\phi]} U_\phi$
onto
$\roman T_{[\phi]}\roman {Rep}_{\xi}(\Gamma,G)$.
Consequently a choice of representative $A$
(in its class $[A]$)
induces an isomorphism of Zariski tangent spaces from
$\roman T_{[0]} \roman H_A$
onto
$\roman T_{[A]} N(\xi)$,
and
a choice of representative $\phi$
(in its class $[\phi]$)
induces an isomorphism of Zariski tangent spaces from
$\roman T_{[0]} \roman H_\phi$
onto
$\roman T_{[\phi]}\roman {Rep}_{\xi}(\Gamma,G)$. \qed
\endproclaim

This reduces the study of the Zariski tangent spaces
to our local models, to which we now turn.
Let $W$ be a finite dimensional unitary
representation of a compact Lie group $K$,
and let $\Theta$ denote its unique momentum mapping
from $W$ to $k^*$ having the value zero at the origin;
further, let $V=\Theta^{-1}(0)$,
with smooth structure $C^{\infty}(V)$ given by (6.1.1)
and
$W_{\roman{red}}=V/K$, its reduced space,
with smooth structure
$C^{\infty}(W_{\roman{red}})$ given by (6.1.2).
By (6.13),
$(W_{\roman{red}},C^{\infty}(W_{\roman{red}}))$
has smooth bump functions arbitrarily close to every point.
Hence, by (7.1),
the inclusion of an arbitrary open connected subset
of
$W_{\roman{red}}$ containing the class $[0]$ of the origin,
with its induced smooth structure,
induces an isomorphism of Zariski tangent spaces
at $[0]$.
The Zariski tangent space
$\roman T_0V$
of $V$ at the origin equals the linear span
$\roman{Vect}(V)$
of
$V$ in $W$, and
projection from
$V$ to
$W_{\roman{red}}$
induces a linear map $\lambda$ from
$\roman T_0V$
to
$\roman T_{[0]}W_{\roman{red}}$.
To deduce information
about the Zariski
tangent space
$\roman T_{[0]}W_{\roman{red}}$,
we denote
the space of $K$-invariants
by $W^K$ and its counter part,
that is the space arising from dividing out the $K$-action,
by
$W_K$.
The kernel $J_K(W)$
of the canonical projection
from $W$ to
$W_K$ is the linear span
of the elements $xw -w,\, x \in K,\, w \in W$.
A little thought reveals that the
orthogonal complement of $W^K$ in $W$
equals the subspace $J_K(W)$, that is,
as a $K$-representation, $W =W^K \oplus J_K(W)$.
Moreover,
the zero locus $V$ contains the subspace $W^K$
of $K$-invariants, and the projection from
$V$ to
$W_{\roman{red}}$, restricted to
$W^K$, is a homeomorphism identifying
the latter with the (smooth)
stratum $S$ in which the class $[0]$ of the origin lies.
The (smooth) tangent space $\roman T_{[0]}S$ of $S$ at
$[0]$ is thus just a copy of
$W^K$, and the inclusion of $S$ into
$W_{\roman{red}}$
induces an injection of
$\roman T_{[0]}S \cong W^K$
into
$\roman T_{[0]}W_{\roman{red}}$.
Furthermore,
with respect to the decomposition
into connected components of
orbit types,
the algebra of invariants $ (C^{\infty}(W))^K$
endows
the orbit space $W\big /K$ with a smooth structure
$C^{\infty}(W\big /K)$,
and the inclusion of
$W_{\roman{red}}$
into
$W\big /K$ is smooth.
Since the induced map
from
$C^{\infty}(W\big /K)$
to
$C^{\infty}(W_{\roman{red}})$
is surjective,
the derivative
$\roman T_{[0]}W_{\roman{red}}
\to
\roman T_{[0]}(W\big /K)$
of this inclusion
is in fact injective.

\proclaim{Lemma 7.4}
Suppose the zero locus $V$
of $\Theta$
spans
$W$
so that the Zariski tangent space
$\roman T_0V$ equals $W$
whence the
linear map
$\lambda$ then goes from
$W$ to
$\roman T_{[0]}W_{\roman{red}}$.
Then $\lambda$
has kernel
$J_K(W)$
and image equal to the
(smooth)
tangent space
$\roman T_{0} S$, viewed as a subspace
of $\roman T_{[0]}W_{\roman{red}}$.
In particular, $\lambda$ is injective and hence an isomorphism
if and only if $W$ is a trivial $K$-representation.
\endproclaim

In the language of
p.~71 of
\cite\armgojen,
the condition says that,
$V$ being viewed as a {\it constraint set\/},
the \lq\lq spanning condition\rq\rq\
is satisfied at  $0\in V$.

\demo{Proof}
View $W$ as a real vector space,
consider
the algebra
$\bold R[W]$
of real polynomials on $W$,
and pick a finite
set of homogeneous
generators
$(\kappa_1,\dots,\kappa_k)$
of the subalgebra
$\bold R[W]^K$
of $K$-invariant polynomials.
Then the Hilbert
map
$\kappa$ from
$W$
to $\bold R^k$
which assigns
$\kappa(w) = (\kappa_1(w),\dots,\kappa_k(w))$
to a vector
$w \in W$
descends to an injective map
$\widetilde \kappa$ from
$W/K$ to $\bold R^k$.
In view of a result
of
\cite\gwschwar,
with reference to the  smooth structure
$C^{\infty}(W/K)$,
the map $\widetilde \kappa$ is proper,
that is, the induced map
from
$C^{\infty}(\bold R^k)$
to
$C^{\infty}(W/K)$ is surjective,
and hence the derivative of
$\widetilde \kappa$
at
the orbit
$[0] = 0\cdot K$
is injective;
further, when the number $k$ is minimal,
by a result of
\cite\mathone,
this derivative
is even an isomorphism from
$\roman T_{[0]}(W/K)$ onto $\bold R^k$.
Thus, for $k$ minimal,
the
canonical map
from $W$ to
$\roman T_{[0]}(W\big /K)$
comes down to the derivative
${
d\kappa(0)
\colon
W
@>>>
\bold R^k
}$
of the Hilbert map at the origin,
and the latter
decomposes
into the linear map $\lambda$
from
$W$ to
$\roman T_{[0]}W_{\roman{red}}$
and
the injection
from
$\roman T_{[0]}W_{\roman{red}}$
into
$\roman T_{[0]}(W\big /K)$
which embeds
$\roman T_{[0]}W_{\roman{red}}$
into a $k$-dimensional
vector space.
However,
$W =W^K \oplus J_K(W)$, and
$d\kappa(0)$ vanishes on
$J_K(W)$ and identifies
$W^K$ with a subspace of $\bold R^k$,
in fact, with what corresponds to the tangent space
$\roman T_{0} S$.
In particular,
$\lambda$  to be injective means that
$W^K$ equals $W$, that is to say, that $K$ acts trivially on $W$. \qed
\enddemo

Next we recall the following well known fact.

\proclaim{Proposition 7.5}
As a smooth space,
$W_{\roman{red}}$
is semi-algebraic.
\endproclaim

We reproduce a proof, for reference in the next Section.

\demo{Proof}
After a choice of invariant polynomials
$(\kappa_1,\dots,\kappa_k)$ has been made,
by the Tarski-Seidenberg theorem,
the resulting injective map
$\widetilde \kappa$
from
$W/K$ to $\bold R^k$
realizes
$W/K$ as a semi-algebraic subset of $\bold R^k$,
in fact, of the real affine categorical
quotient
$W//K$, that is,
of the real affine variety
determined by a finite set of relations for
the algebra of invariants $\bold R[W]^K$.
The composite of
$\widetilde \kappa$
with the canonical injection
of $W_{\roman{red}}$ into $W/K$
embeds
$W_{\roman{red}}$
into $\bold R^k$.
To see this embedding is semi-algebraic,
write $I_V$
for the ideal of $V$ in $\bold R[W]$ and
consider  the real affine
coordinate ring $A[V]= \bold R[W]/I_V $ of $V$.
Since $K$ is compact,
the canonical map
from
$\bold R[W]^K\big / I_V^K$
to the $K$-invariants
$A[V]^K$ is an isomorphism.
Let $\phi_1,\dots,\phi_{\ell}$
be a finite set of generators
of $I_V^K$;
when we write them out in the generators
$(\kappa_1,\dots,\kappa_k)$,
we obtain a polynomial map
$\Phi$ from $\bold R^k$ to $\bold R^{\ell}$
so that
$W_{\roman{red}}$
amounts to the intersection of
$W/K$
with the real affine set
$\Phi^{-1}(0)$
whence
$W_{\roman{red}}$
is semi-algebraic in $\bold R^k$. \qed
\enddemo
\smallskip\noindent
{\smc Remark.}
We have seen above that the
inclusion of
$W_{\roman{red}}$ into
$W/K$ induces an embedding
of the Zariski tangent space
$\roman T_{[0]}W_{\roman{red}}$
into the Zariski tangent space $\roman T_{[0]}(W/K)$.
The above embedding of
$W_{\roman{red}}$
into $\roman T_{[0]}(W/K)$
passes to an embedding into
$\roman T_{[0]}W_{\roman{red}}$.
In fact, the embedding of
$W/K$ into its Zariski tangent space
is induced by the canonical embedding of
$W$ into its tangent space
$\roman T_{0}W$
which assigns to a vector $w \in W$ its directional derivative
at the origin
on smooth functions on $W$.
It is obvious that this association passes
to one which assigns
to a vector $w \in V$
an element in the Zariski tangent space
$\roman T_{0}V$, viewed as a linear subspace of
$\roman T_{0}W$,
and hence, by $K$-invariance,
to an embedding
of $W_{\roman{red}}$
into
its Zariski tangent space
$\roman T_{[0]}W_{\roman{red}}$
as a semi-algebraic set.
An example will be examined in the next Section.

\smallskip
We now
apply the above to moduli spaces.
For a central Yang-Mills connection $A$,
we shall denote by $V_A$
the zero locus of
the quadratic mapping
$\Theta_A$ from
$\roman H_A^1 $ to $\roman H_A^2$,
cf. Section 6, and likewise,
for
$\phi$
in $\roman {Hom}_{\xi}(\Gamma,G)$,
by
$V_\phi$
the zero locus of
the quadratic mapping
$\Theta_\phi$ from
$\roman H_\phi^1 $ to $\roman H_\phi^2$.

\proclaim{Lemma 7.6}
For every central Yang-Mills connection $A$,
the zero locus
$V_A$
spans $\roman H_A^1(\Sigma,\roman{ad}(\xi))$.
Likewise,
for every
$\phi$
in $\roman {Hom}_{\xi}(\Gamma,G)$,
the zero locus
$V_\phi$
spans $\roman H^1(\pi,g_\phi)$.
\endproclaim

The proof of this Lemma
requires some preparation.
We shall denote by
$\Cal N(\xi)^-$
the subspace
of
central Yang-Mills connections $A$
having the property
that the Lie bracket $[\cdot,\cdot]_A$
is zero on $\roman H^1_A$.
Recall that
a description the space $\Cal A_{A}(\xi)$
for a central Yang-Mills connection $A$
has been reproduced in Section 6 above.
It is proved in \cite\singula~(2.8)
that, near a central Yang-Mills connection $A$, the space $\Cal N(\xi)$
coincides with
$\Cal A_{A}(\xi)$
and hence
is smooth
near $A$,
with tangent space
$\roman T_{A} \Cal N(\xi)$
equal to the space $Z_{A}^1(\Sigma,\roman{ad}(\xi))$
of 1-cocycles
if and only if
$A$ lies in $\Cal N(\xi)^-$.
Thus the subspace $\Cal N(\xi)^-$
is a smooth submanifold of $\Cal A(\xi)$,
and
from the exactness
of (6.4)
we deduce that,
for every point $A$ of $\Cal N(\xi)^-$,
the operator of covariant derivative
$d_{A}$ gives rise to the exact sequence
$$
0
@>>>
\roman T_{A} \Cal N(\xi)
@>>>
\roman T_{A} \Cal A(\xi)
@>{d_{A}}>>
\Omega^2(\Sigma,\roman{ad}(\xi))
@>>>
\roman H_{A}^2(\Sigma,\roman{ad}(\xi))
@>>>
0
\tag7.7
$$
of real vector spaces.
In fact,
the points of $\Cal N(\xi)^-$
are exactly the weakly regular points
(p. 300 of \cite\abramars)
for the momentum mapping $J$ from
$\Cal A(\xi)$ to
$\Omega^2(\Sigma,\roman{ad}(\xi))$,
cf. Section 6.
\smallskip
Denote by
$\Cal N^{\roman{top}}(\xi)$
the subspace of
$\Cal N(\xi)$
which consists
of central Yang-Mills connections
$A$
having the property that
$Z_A$
acts trivially on
$\roman H_A^1(\Sigma,\roman{ad}(\xi))$, so that
the top stratum
$N^{\roman{top}}(\xi)$
equals $\Cal N^{\roman{top}}(\xi)\big/\Cal G(\xi)$,
see our paper~\cite\singulat.
By  \cite\singulat~(1.5),
there is a certain subgroup $Z^{\roman{top}}$ of $G$,
unique up to conjugacy,
such that
under (1.1)
the image
of the stabilizer
$Z_A$
of every
central Yang-Mills connection
$A$
in $\Cal N^{\roman{top}}(\xi)$
is conjugate to
$Z^{\roman{top}}$.
Since $\Theta_A$ is a momentum mapping
for every central Yang-Mills connection $A$,
$\Cal N^{\roman{top}}(\xi)$
is a subspace of $\Cal N^-(\xi)$,
in fact,
a smooth codimension zero submanifold
since for every
$A \in \Cal N^{\roman{top}}(\xi)$
the tangent map of
the inclusion
$\Cal N^{\roman{top}}(\xi) \subseteq \Cal N^-(\xi)$
amounts to  the identity mapping
of $Z_A^1(\Sigma,\roman{ad}(\xi))$.
\smallskip
In what follows, by the {\it dimension\/}
$\dim V_A$ of $V_A$ we mean the dimension
of its non-singular part
$V^-_A \subseteq V_A$.

\demo{Proof of {\rm (7.6)}}
Since the top stratum $N^{\roman{top}}(\xi)$
is dense in $N(\xi)$,
cf.
\cite\singulat\ (1.4),
arbitrarily close to $[A]$ there is
a point $[\widetilde A]$
in the top stratum,
and
we may assume that the group $Z^{\roman{top}}$ is the stabilizer
$Z_{\widetilde A}$
of $\widetilde A$.
Then a neighborhood of
the point $x$ of $V_A$ corresponding to
$\widetilde A$ is the total space of a
$Z_A$-fibre bundle,
having as base space a neighborhood of
the class $[x]$
in $V_A/Z_A$
and as fibre the homogeneous space
$Z_A/Z^{\roman{top}}$.
Consequently
$$
\align
\dim V_A
=
\dim \roman T_x V_A
&
= \dim N^{\roman{top}}(\xi)
+
\dim Z_A - \dim Z^{\roman{top}}
\\
&
= \dim \roman H^1_{\widetilde A}+
\dim Z_A - \dim Z^{\roman{top}} .
\endalign
$$
However,
for every central Yang-Mills connection $\overline A$,
the
twisted integration mapping yields an isomorphism from
$\roman H_{\overline A}^*(\Sigma,\roman{ad}(\xi))$
onto $\roman H^*(\pi,g_{\rho\overline A})$.
Now
an Euler characteristic argument
in the chain complex calculating
the corresponding group cohomologies
establishes
equality between the two alternating sums
$\dim\roman H_A^0
-\dim\roman H_A^1
+\dim\roman H_A^2
$
and
$\dim\roman H_{\widetilde A}^0
-\dim\roman H_{\widetilde A}^1
+\dim\roman H_{\widetilde A}^2$.
Since
$\dim \roman H_A^2
=
\dim \roman H_A^0=
\dim Z_A
$
and
${\dim \roman H_{\widetilde A}^2
=
\dim \roman H_{\widetilde A}^0=
\dim Z^{\roman{top}}}$,
we conclude
$$
\dim\roman H_{\widetilde A}^1
-
2 \dim Z^{\roman{top}}
=
\dim\roman H_{A}^1
-
2\dim Z_A,
\tag7.6.1
$$
and thence
$$
\dim V_A
=
\dim \roman H_A^1
-
\dim Z_A + \dim Z^{\roman{top}}.
\tag7.6.2
$$
\smallskip

Next we assert that,
at the image $x$ of $[\widetilde A]$
in $V_A\subseteq \roman H^1_A$,
the derivative
${
d\Theta_A(x) \colon
\roman H_A^1
\to \roman H_A^2
}$
of $\Theta_A$
has rank
$$
\roman{rank}
\left(d\Theta_A(x)\right) =
\dim
\roman H_A^2 - \dim Z^{\roman{top}}
=
\dim Z_A- \dim Z^{\roman{top}}.
\tag7.6.3
$$
Now, at a point
$\overline A \in \Cal A_A$,
the smooth submanifold
$
\Cal A_A
$
of $\Cal A(\xi)$
has tangent space
$$
\roman T_{\overline A} \Cal A_A
=
\{ \phi; d_{\overline A} \phi \in \Cal H_A^2(\Sigma,\roman{ad}(\xi))\}
\subseteq
\Omega^1(\Sigma,\roman{ad}(\xi)).
$$
In other words,
the right-hand unlabelled arrow being the inclusion,
the square
$$
\CD
\roman T_{\overline A} \Cal A_A
@>{d_{\overline A}|}>>
\Cal H_A^2(\Sigma,\roman{ad}(\xi))
\\
@VVV
@VVV
\\
\roman T_{\overline A} \Cal A(\xi)
@>{d_{\overline A}}>>
\Omega^2(\Sigma,\roman{ad}(\xi))
\endCD
\tag7.6.4
$$
is a pull back diagram.
By construction,
${
\Cal N(\xi)
=
\{\widehat A \in \Cal A_A; K_{\widehat A } = K_{\xi} \};
}$
here
$K_{\xi}$
refers to the element
of $\Cal H_A^2(\Sigma,\roman{ad}(\xi))$
determined
by the topology of $\xi$,
see Section 2 of \cite\singula.
Since (7.6.4) is a pull back diagram, by standard principles,
at a
point
$\widehat A$ of $\Cal N^-(\xi)$
the
sequence (7.7)
induces an exact sequence of real vector spaces
$$
0
@>>>
\roman T_{\widehat A} \Cal N(\xi)
@>>>
\roman T_{\widehat A} \Cal A_A
@>{d_{\widehat A}}>>
\Cal H_A^2(\Sigma,\roman{ad}(\xi))
@>>>
\roman H_{\widehat A}^2(\Sigma,\roman{ad}(\xi)).
\tag7.6.5
$$
Notice at present we cannot assert that
the last arrow in (7.6.5) is surjective.
\smallskip

Next we recall that,
for
$\widehat A$
in $\Cal N^-(\xi)$
and close to $A$,
the smooth submanifold
$\Cal M_A$ of
$\Cal A_A$,
cf. \cite\singula~(2.16) and Section 6 above,
has tangent space
$
\roman T_{\widehat A} \Cal M_A
$
equal to
$\roman T_{\widehat A} \Cal A_A
\cap
\roman{ker}(d_A^*)$;
hence such a point
$\widehat A$ gives rise to
the exact sequence
$$
0
@>>>
\left(\roman T_{\widehat A} \Cal N(\xi)
\cap
\roman{ker}(d_A^*)\right)
@>>>
\roman T_{\widehat A} \Cal M_A
@>>>
d_{\widehat A}\left(\roman T_{\widehat A} \Cal M_A\right)
@>>>
0
$$
which,
cf. Section 2 of \cite\singula,
with
$\Cal N_A = \Cal N(\xi) \cap \Cal M_A$,
looks like
$$
0
@>>>
\roman T_{\widehat A} \Cal N_A
@>>>
\roman T_{\widehat A} \Cal M_A
@>>>
d_{\widehat A}\left(\roman T_{\widehat A} \Cal M_A\right)
@>>>
0.
\tag7.6.6
$$
We note that, near $A$,
$\Cal N_A$ also  equals the intersection
$\Cal N(\xi) \cap (A+\roman{ker}(d_A^*))$.
\smallskip

Let now $\widetilde A$ be a point
close to $A$
representing a point of $N^{\roman{top}}$;
then
$\widetilde A$
lies in particular in $\Cal N^-(\xi)$,
and
near $\widetilde A$,
the restriction to
$\Cal N_A$
of the projection map from
$\Cal N(\xi)$ onto $N(\xi)$ is
a fibre bundle map onto its image,
having fibre the homogeneous space
$Z_A\big/Z_{\widetilde A}$.
Consequently, in view of (7.6.1),
$$
\align
\dim \Cal N_A
&=
\dim N(\xi) +
\dim Z_A- \dim Z_{\widetilde A}
\\
&=
\dim \roman H_{\widetilde A}^1 +
\dim Z_A- \dim Z_{\widetilde A}
\\
&=
\dim \roman H_{A}^1 +
\dim Z_{\widetilde A} -\dim Z_A .
\endalign
$$
However,
$\dim \Cal M_A
=
\dim \roman H_A^1$.
Consequently
$$
\align
\dim
d_{\widetilde A}\left(\roman T_{\widetilde A} \Cal M_A\right)
&=
\dim \roman H_A^1
-
\dim \Cal N_A
\\
&=
\dim Z_A- \dim Z_{\widetilde A}
\\
&=
\dim \roman H_A^2
-
\dim \roman H_{\widetilde A}^2,
\endalign
$$
whence the
exact sequence (7.6.5) furnishes
the exact sequence
$$
0
@>>>
\roman T_{\widetilde A} \Cal N_A
@>>>
\roman T_{\widetilde A} \Cal M_A
@>{d_{\widetilde A}}>>
\Cal H_A^2(\Sigma,\roman{ad}(\xi))
@>>>
\roman H_{\widetilde A}^2(\Sigma,\roman{ad}(\xi))
@>>>
0
\tag7.6.7
$$
of finite dimensional real vector spaces;
notice its exactness at
$\roman T_{\widetilde A} \Cal M_A$
is implied by that of (7.6.6).
By construction,
the Kuranishi map identifies (7.6.7) with the sequence
$$
0
@>>>
\roman T_{x} V_A
@>>>
\roman T_{x} \roman H_A^1(\Sigma,\roman{ad}(\xi))
@>{d\Theta_{A}(x)}>>
\roman H_A^2(\Sigma,\roman{ad}(\xi))
@>>>
\roman H_{\widetilde A}^2(\Sigma,\roman{ad}(\xi))
@>>>
0
\tag7.6.8
$$
which is therefore exact.
In particular, the point $x\in V_A$ is  weakly regular
for
$\Theta_A$,
and hence
$d\Theta_{A}(x)$
has rank asserted in (7.6.3).
\smallskip

Finally we show that the latter implies that the
{\it real\/} linear span
$\roman{Vect}(V_A)$
of $V_A$ in $\roman H_A^1(\Sigma,\roman{ad}(\xi))$
equals
the whole space
$\roman H_A^1(\Sigma,\roman{ad}(\xi))$.
In fact,
the cone $V_A$ is obviously stable under
$Z_A$.
Moreover,
in view of
\cite\singula\ (2.27),
for every $\eta \in \Cal H^1_A(\Sigma,\roman{ad}(\xi))$,
the value $[\eta,\eta] \in \Cal H^2_A(\Sigma,\roman{ad}(\xi))$
is zero if and only if
$[*\eta,*\eta]= 0$;
here $*$ refers to the corresponding duality operator,
cf.
\cite\singula\ (1.1.5).
Consequently the cone $V_A$ is  stable under
the duality operator $*$.
However this duality operator
induces the {\it
complex\/} structure
on
$\roman H_A^1(\Sigma,\roman{ad}(\xi))$.
Hence the {\it real\/} linear span
$\roman{Vect}(V_A)$
of $V_A$ in $\roman H_A^1(\Sigma,\roman{ad}(\xi))$
equals
its {\it complex\/} linear span
in $\roman H_A^1(\Sigma,\roman{ad}(\xi))$;
the complex vector space $\roman{Vect}(V_A)$
thus inherits a structure of a
unitary
$Z_A$-representation,
and as a unitary
$Z_A$-representation,
the space
$\roman H_A^1(\Sigma,\roman{ad}(\xi))$
decomposes
into the direct sum
of
$\roman{Vect}(V_A)$
and its orthogonal complement
$\roman{Vect}(V_A)^{\bot}$.
Moreover the restrictions
$\Theta^1_A$
and
$\Theta^2_A$
of
$\Theta_A$
to
$\roman{Vect}(V_A)$
and
$\roman{Vect}(V_A)^{\bot}$, respectively,
are the unique momentum mappings
for these unitary
$Z_A$-representations
having the value zero
at the origin.
By construction, the cone
$V_A$ lies in the summand
$\roman{Vect}(V_A)$,
whence the zero locus
${
(\Theta^2_A)^{-1}(0) \subseteq \roman{Vect}(V_A)^{\bot}
}$
consists merely of the origin.
Hence, whatever
weakly regular point
$x$ of $V_A$,
the rank of the derivative
$d\Theta_A (x)$
coincides with the rank of the restriction
$
d\Theta^1_A (x)
$
to
$\roman T_x \roman{Vect}(V_A)= \roman{Vect}(V_A)$.
Consequently
${
\dim V_A
=
\dim\roman{Vect}(V_A)
-
\dim Z_A+ \dim  Z^{\roman{top}}.
}$
However,
in view of (7.6.2),
this can only happen if
${
\dim\roman{Vect}(V_A)
=
\dim
\roman H_A^1(\Sigma,\roman{ad}(\xi)),
}$
whence
\linebreak
${
\roman{Vect}(V_A)
=
\roman H_A^1(\Sigma,\roman{ad}(\xi))
}$
as asserted. \qed
\enddemo

\smallskip
\noindent
{\smc Remark 7.8.}
Let $K$ be a compact Lie group, with Lie algebra $k$,
let $W$ be an
$n$-dimensional
unitary representation of $K$,
and let
$\mu$
be the unique momentum mapping
from $W$ to  $k^*$
having the value zero at the origin.
Its derivative
at the origin
is zero,
the kernel of
$d\mu(0)$ in fact equals the whole space $W$,
and
the Zariski tangent space
$\roman T_0(\mu^{-1}(0))$
at the origin
of the zero locus
$\mu^{-1}(0)$
is obviously a subspace
of
the kernel of
$d\mu(0)$.
However
in general
the Zariski tangent space
does {\it not\/} coincide
with the kernel of
$d\mu(0)$.
To see this,
suppose that the irreducible representations
in $W$ are all non-trivial, that
$K$ is a subgroup of the unitary group $U(n)$,
and that $K$ contains the central circle subgroup
$S^1$ of $U(n)$.
Since the momentum mapping
for the $S^1$-action on $\bold C^n$
is given by the assignment to
$\bold z \in \bold C^n$ of $||\bold z||^2$,
the zero level set
$\mu^{-1}(0)$
will then consist of the origin only,
the Zariski tangent space
of which is of course trivial.
Thus (7.6)
is {\it non-trivial\/}.
\smallskip

The decompositions
of $N(\xi)$ and $\roman {Rep}_{\xi}(\Gamma,G)$
into connected components of orbit types
have been shown to be a stratification
in \cite\singulat.
If $[A]$ lies in the stratum
$N_{(K)}$,
the inclusion of
$N_{(K)}$ into  $N(\xi)$
induces an injection
$\roman T_{[A]}(N_{(K)}) \to \roman T_{[A]}N(\xi)$
of Zariski tangent spaces,
and
$\roman T_{[A]}(N_{(K)})$ will in this way be viewed
as a linear subspace of
$\roman T_{[A]}N(\xi)$;
this is e.~g. a consequence of (6.2) combined with (7.4).
Notice
$\roman T_{[A]}N_{(K)}$
amounts to the usual smooth tangent space
of the smooth manifold
$N_{(K)}$.
It is clear that the same kind of remarks
can be made
for an arbitrary point $[\phi]$ of
$\roman {Rep}_{\xi}(\Gamma,G)$
and the stratum
$\roman {Rep}_{\xi}(\Gamma,G)_{(K)}$
in which it lies.
A point in the top stratum
$N^{\roman{top}}(\xi)$
will be referred to as a {\it non-singular\/} point
of
$N(\xi)$, cf. \cite\singulat.
Accordingly
the representation space
$\roman{Rep}_{\xi}(\Gamma,G)$
has a
{\it non-singular\/} part or
{\it top\/} stratum
$\roman{Rep}^{\roman{top}}_{\xi}(\Gamma,G)$,
and a point in
$\roman{Rep}^{\roman{top}}_{\xi}(\Gamma,G)$
will be said to be a {\it non-singular\/} point
of
$\roman{Rep}_{\xi}(\Gamma,G)$.
We now collect a number of consequences
of the above results.
\smallskip\noindent
{\bf 7.9.}
Let
$[A]$ be a point of $N(\xi)$.
In view of (6.2) and (7.6), a choice of representative $A$
in its class $[A]$ determines a linear map
$\lambda_A$ from
$\roman H_A^1(\Sigma,\roman{ad}(\xi))$
to
$\roman T_{[A]}N(\xi)$.
In fact,
this map is the composite of the linear map
$\lambda$ from
$\roman H_A^1(\Sigma,\roman{ad}(\xi))$
to
$\roman T_{[0]} W_A$,
cf. (7.4),
with the derivative of the injection
of $W_A$ into $N(\xi)$ given in (6.2),
where
$W_A$
refers to an open neighborhood of the class of zero in
$\roman H_A$ of the kind coming into play in (6.2).
By construction,
$\lambda_A$
depends on the Kuranishi map;
however the latter, in turn,
depends merely on the data
coming into play in the definition of
$N(\xi)$.
It is in this sense that
a choice of representative
$A$
of $[A]$
in fact {\it determines\/}
$\lambda_A$.
The map $\lambda_A$ has the following properties:
\newline\noindent {\rm (1)}
{\sl It is independent of the choice of
$A$ in the sense that,
for every gauge transformation
$\gamma \in \Cal G(\xi)$,
the composite
$$
\roman H_{A}^1(\Sigma,\roman{ad}(\xi))
@>{\gamma_{\sharp}}>>
\roman H_{\gamma (A)}^1(\Sigma,\roman{ad}(\xi))
@>{\lambda_{\gamma (A)}}>>
\roman T_{[A]}N(\xi)
$$
of the induced linear isomorphism $\gamma_{\sharp}$
with $\lambda_{\gamma (A)}$
coincides with
$\lambda_A$.
\newline\noindent {\rm (2)}
Its kernel
equals the subspace
$J_K(\roman H_A^1)$
of $\roman H_A^1 =\roman H_A^1(\Sigma,\roman{ad}(\xi))$,
where
$K = Z_A$, the stabilizer of $A$.
\newline\noindent {\rm (3)}
Its image equals the
(smooth)
tangent space
$\roman T_{[A]}(N_{(K)})$, viewed as a subspace
of
$\roman T_{[A]}N(\xi)$
in a sense explained above,
where
$N_{(K)}$ denotes the stratum in which $[A]$ lies.
\newline\noindent {\rm (4)}
It is an isomorphism
if and only if $[A]$ is a non-singular point of $N(\xi)$.}
\smallskip
These follow at once from
(7.4) except statement (1)
the proof of which we leave to the reader.

\smallskip\noindent
{\bf 7.10.}
Let
$[\phi]$
be a point of $\roman {Rep}_{\xi}(\Gamma,G)$.
In view of (6.3) and (7.6), a choice of representative $\phi$
in $\roman {Hom}_{\xi}(\Gamma,G)$
in its class $[\phi]$ determines a linear map
$\lambda_\phi$
from
$\roman H^1(\pi,g_{\phi})$
to
$\roman T_{[\phi]} \roman {Rep}_{\xi}(\Gamma,G)$.
In fact,
this map is the composite of the linear map
$\lambda$ from
$\roman H^1(\pi,g_{\phi})$
to
$\roman T_{[0]} W_\phi$,
cf. (7.4), with the derivative of the injection
of $W_\phi$ into $\roman {Rep}_{\xi}(\Gamma,G)$ given in (6.3),
where
$W_\phi$
refers to an open neighborhood of the class of zero in
$\roman H_\phi$ of the kind coming into play in (6.3).
By construction,
the injection
of $W_\phi$ into $\roman {Rep}_{\xi}(\Gamma,G)$
depends a priori on the Kuranishi map and in particular
on the choice of Riemannian metric on $\Sigma$.
However
$\lambda_\phi$
does {\it not\/}
depend on this choice.
In fact, by (5.7),
the derivative of the Wilson loop mapping
$\rho$ from $\Cal A(\xi)$ to $\roman {Hom}(F,G)$
at a central Yang-Mills connection $A$, restricted to
the subspace $Z^1_A(\Sigma,\roman{ad}(\xi))$
of 1-cocycles in
$\Omega^1(\Sigma,\roman{ad}(\xi))= \roman T_A\Cal A(\xi)$,
amounts to the composite
$$
Z^1_A(\Sigma,\roman{ad}(\xi))
@>{\roman{Int}_A|}>>
Z^1(\pi,g_{\phi})
@>{\roman L_{\phi}}>>
\roman T_{\phi}\roman{Hom}_{\xi}(\Gamma,G)
$$
of the restriction
$\roman{Int}_A|$ of the twisted
integration mapping
with left translation
$\roman L_{\phi}$
from
$Z^1(\pi,g_{\phi})$
to
$\roman T_{\phi}\roman{Hom}_{\xi}(\Gamma,G)$,
whatever Riemannian metric on $\Sigma$;
here $\phi =\rho(A) \in \roman{Hom}_{\xi}(\Gamma,G)$.
Since every
$\phi \in \roman{Hom}_{\xi}(\Gamma,G)$
arises in this way,
for every such $\phi$,
the diagram
$$
\CD
Z^1(\pi,g_{\phi})
@>{\roman L_{\phi}}>>
\roman T_{\phi}\roman{Hom}_{\xi}(\Gamma,G)
\\
@VVV
@VVV
\\
\roman H^1(\pi,g_{\phi})
@>{\lambda_{\phi}}>>
\roman T_{\phi}\roman{Rep}_{\xi}(\Gamma,G)
\endCD
$$
is commutative, the unlabelled
vertical maps being the obvious ones.
Hence
a choice of representative $\phi$
in its class $[\phi]$ indeed determines a linear map
$\lambda_\phi$
as asserted which does {\it not\/} depend on
a choice of Riemannian metric on $\Sigma$.
The map $\lambda_\phi$ has the following properties:
\newline\noindent {\rm (1)}
{\sl It is independent of the choice of
$\phi$ in the sense that,
for every $x\in G$, the
composite
$$
\roman H^1(\pi,g_{\phi})
@>{\roman{Ad}_\flat(x)}>>
\roman H^1(\pi,g_{x\phi})
@>{\lambda_{x\phi}}>>
\roman T_{[\phi]} \roman {Rep}_{\xi}(\Gamma,G)
$$
of the induced linear isomorphism $\roman{Ad}_\flat(x)$
with
$\lambda_{x\phi}$
coincides with
$\lambda_{\phi}$.
\newline\noindent {\rm (2)}
Its kernel equals the subspace
$J_K(\roman H^1(\pi,g_{\phi}))$
of $\roman H^1(\pi,g_{\phi})$,
where
$K = Z_{\phi}$, the stabilizer of ${\phi}$.
\newline\noindent {\rm (3)}
Its image equals the
(smooth)
tangent space
$\roman T_{[{\phi}]}(\roman {Rep}_{\xi}(\Gamma,G)_{(K)})$, viewed as a subspace
of
$\roman T_{[{\phi}]}\roman {Rep}_{\xi}(\Gamma,G)$
in a sense explained above,
where
$\roman {Rep}_{\xi}(\Gamma,G)_{(K)}$ denotes the stratum in which $[\phi]$
lies.
\newline\noindent {\rm (4)}
It is an isomorphism
if and only if $[{\phi}]$ is a non-singular point of
$\roman {Rep}_{\xi}(\Gamma,G)$.}
\smallskip
These follow again at once from
(7.4) except statement (1)
the proof of which is formally the same as that of (7.9(1)).
\smallskip
The statements of (7.9) and (7.10) are related by the
fact that, for every central Yang-Mills connection $A$,
the diagram
$$
\CD
\roman H_A^1(\Sigma,\roman{ad}(\xi))
@>{\lambda_A}>>
\roman T_{[A]} N({\xi})
\\
@V{\roman{Int}_A}VV
@V{d\rho_{\flat}[A]}VV
\\
\roman H^1(\pi,g_{\rho(A)})
@>>{\lambda_{\rho(A)}}>
\roman T_{[\rho(A)]} \roman {Rep}_{\xi}(\Gamma,G)
\endCD
\tag7.11
$$
is commutative.
Thus at a non-singular point $[A]$ of $N(\xi)$, the derivative of the
Wilson loop mapping comes down to the twisted integration
mapping $\roman{Int}_A$ from
$\roman H_A^1(\Sigma,\roman{ad}(\xi))$
to
$\roman H^1(\pi,g_{\rho(A)})$.

\smallskip
\noindent
{\smc Remark 7.12.}
At a singular point
$[\phi]$ of
$\roman{Rep}_{\xi}(\Gamma,G)$,
the Zariski tangent space
$\roman T_{[\phi]} \roman{Rep}_{\xi}(\Gamma,G)$
with respect to the smooth structure
$C^{\infty}(\roman{Rep}_{\xi}(\Gamma,G))$
does {\it not\/}
boil down to
$\roman H^1(\pi,g_{\phi})$,
cf. what is said on p.~205 of \cite\goldmone.
An example where this phenomenon really occurs
will be given in the next Section.
\smallskip

Here is an immediate consequence of (6.2), (6.3), and (7.5):

\proclaim{Theorem 7.13}
As smooth spaces, $N(\xi)$
and its diffeomorphe $\roman{Rep}_{\xi}(\Gamma,G)$
are locally semi-algebraic.
\endproclaim

Next we spell out our {\it fifth main result\/}.
For every $\phi \in \roman{Hom}_{\xi}(\Gamma,G)$,
the kernel of the derivative
$dr_{\phi}$
from
$\roman T_{\phi}\roman{Hom}(F,G)$ to
$\roman T_{\roman{exp}(X_{\xi})}G$,
with reference to the word map
$r$ from $\roman{Hom}(F,G)$ to $G$,
yields a notion of
{\it not necessarily reduced\/}
Zariski tangent space,
and it is clear that
the Zariski
tangent space $\roman T_{\phi}\roman{Hom}_{\xi}(\Gamma,G)$
with reference to
the smooth structure $C^{\infty}(\roman{Hom}_{\xi}(\Gamma,G))$
(introduced in Section 3)
is a subspace thereof;
however, a priori
the two spaces
should
{\it not\/}
be confused.

\proclaim{Theorem 7.14}
For every point $\phi \in \roman{Hom}_{\xi}(\Gamma,G)$,
the Zariski tangent space
with reference to
$C^{\infty}\left(\roman{Hom}_{\xi}(\Gamma,G)\right)$
coincides with the
kernel of
the derivative
$dr_{\phi}$.
\endproclaim

Thus our {\it reduced\/} Zariski tangent space
coincides with the other notion of Zariski tangent space.
However, we do not know
whether the ideal
in $C^{\infty}(\roman{Hom}(F,G))$
corresponding to the word map $r$
coincides with its real radical.

\demo{Proof}
Let $A$ be a central Yang-Mills connection
so that $\rho(A) = \phi$.
The smooth  $G$-invariant immersion (6.11.3)
identifies the subspace
${
G \times_{Z_A} \vartheta_A^{-1}(0)
}$
with a $G$-invariant neighborhood
of $\phi$ in
$\roman{Hom}_{\xi}(\Gamma,G)$;
here $\vartheta_A$ refers to the momentum mapping coming into play in Section 6
above.
However the Kuranishi map
$\Phi_A$, cf. \cite\singula\  (2.29) and what is said in Section 6 above,
identifies the inclusion
of
$
G \times_{Z_A} \vartheta_A^{-1}(0)
$ into
$G \times_{Z_A} \Cal M_A$
with the inclusion
of
$G \times_{Z_A}  V_A$
into
$G \times_{Z_A} \roman H^1_A$
where
$V_A \subseteq \roman H^1_A$
refers to the cone $\Theta_A^{-1}(0)$;
see Section 6 above
for any unexplained notation.
Now the tangent space
$\roman T_{[e,0]}
\left(G \times_{Z_A} \roman H^1_A\right)
$
decomposes into a direct sum of
$B^1(\pi,g_{\phi})$ and
$\roman H^1_A$,
that is, it
amounts to the space
$Z^1(\pi,g_{\phi})$ of 1-cocycles,
and in suitable coordinates
near the point
$[e,0]$,
the space
$G \times_{Z_A}  V_A$
boils down  to the zero locus of the
composition of the projection
from $Z^1(\pi,g_{\phi})$
to
$\roman H^1(\pi,g_{\phi})$
with the momentum mapping
$\Theta_A$  from
$\roman H^1_A$
to $\roman H^2_A$.
In view of
(7.6) above, this implies that the Zariski tangent space
of
$G \times_{Z_A}  V_A$
at the point $[e,0]$
equals the
tangent space
${
\roman T_{[e,0]}
\left(G \times_{Z_A} \roman H^1_A\right)
}$
whence the assertion.  \qed
\enddemo

\medskip\noindent{\bf 8. An example}
\smallskip\noindent
Consider the moduli space $N$ of flat
$\roman {SU}(2)$-connections
for a surface $\Sigma$ of genus 2.
This example is already sufficiently
general
to visualize
the global picture which emerges.
As a space, $N$ is just complex projective 3-space,
by a result of
{\smc Narasimhan-Ramanan}~\cite\naramntw.
However we shall see that, as a smooth space,
with smooth structure (3.6),
it looks rather different.
\smallskip
Write  $G =\roman {SU}(2)$,
and let $Z=\{\pm 1\}$
denote the centre of $G$ and
$T=S^1 \subseteq G$
the standard circle subgroup inside
$G$; it is a maximal torus.
The decomposition
of $N$
according to orbit types of flat connections
has the three pieces $N_{G}$, $N_{(T)}$, and  $N_{Z}$,
where the subscript refers to the conjugacy class
of stabilizer;
we recall that
$N_{G}$ consists of 16 isolated points and  that $N_{(T)}$ is connected.
\smallskip
In view of what is said in Section 6 of our paper \cite\locpois,
near a point of the middle stratum
$N_{(T)}$,
as a smooth space,
$N$ looks like a product of
a standard $\bold R^4$
with
a copy of the
cone
$C = \{(u,v,r); u^2 + v^2 = r^2; r \geq 0\}$,
with smooth structure
induced by the embedding of $C$ in 3-space
with coordinates $(u,v,r)$.
In fact,
the latter arises
as the reduced space
for the diagonal
$\roman{SO}(2,\bold R)$-action
on $W = \bold R^2 \times \bold R^2$ with its obvious symplectic
structure, in the following way:
Let $K= \roman{SO}(2,\bold R)$, and  write elements of $W$ in the form
$w=(q,p) \in \bold R^2 \times \bold R^2$.
The algebra
$\bold R[W]^K$ of invariants is generated by
$qq, \ pp,\  qp,\  |q\,p|$,
and the momentum mapping $\mu$ is given by
$\mu(q,p) = |q\,p|$.
However $\mu$
generates the ideal $I_V$
of polynomials
in $\bold R[W]$
vanishing on the zero locus
$V = \mu^{-1}(0)$;
since $\mu$ is $K$-invariant it also
generates the ideal $I_V^K$
of $K$-invariant polynomials
vanishing on $V$.
Thus
the coordinate ring $A[V]$
has four generators
while
the subalgebra of $K$-invariants
$A[V]^K$ is
generated by
$u = qq -pp, v =2qp, r = qq + pp$,
subject to the relation
$r^2 = u^2 + v^2$.
The real affine categorical quotient $W//K$
is the double cone given by this equation
while the reduced space
$W_{\roman{red}}$
amounts to the positive cone $C$,
with the cone point included.
Moreover it is manifest that
the Zariski tangent space $\roman T_{0}C$ at the cone point 0 has dimension 3.
In fact the invariants $u,v,r$ induce a map
$\lambda$ from  $\bold R^4$ to $\bold R^3$
passing through
a map
$\widetilde\lambda$ from $\bold R^4/K$ to $\bold R^3$;
now $\lambda$ has derivative zero at the origin
while
the derivative of
$\widetilde \lambda$
induces an isomorphism
from $\roman T_{[0]}$ onto $\bold R^3$.
Hence
for
a point $[A]$ of the middle stratum
$N_{(T)}$,
the Zariski tangent space
$\roman T_{[A]} N$ has dimension 4 + 3 = 7.
On the other hand,
the dimension of $\roman H_A^1(\Sigma,\roman{ad}(\xi))$
equals 8, and
the linear map
$\lambda_A$ from
$\roman H_A^1(\Sigma,\roman{ad}(\xi))$
to
$\roman T_{[A]} N$
has rank four since the derivative of $\lambda$
at the origin
has rank zero.
Thus
the Zariski tangent space
$\roman T_{[A]} N$
can
in no way be identified with
the cohomology group
$\roman H_A^1(\Sigma,\roman{ad}(\xi))$,
cf. (7.12) above.
\smallskip
Likewise,
in view of what is said in Section 7 of our paper \cite\locpois,
near any
of the 16 points of
$N_{G}$,
as a smooth space,
$N$ looks like the reduced space for the momentum mapping
$\mu$  from $W=(\bold R^3)^4$
to
the dual of $\roman{so}(3,\bold R)$,
for the diagonal $\roman{SO}(3,\bold R)$-action
on $W$ with its obvious symplectic structure,
the $\roman{SO}(3,\bold R)$-action on $\bold R^3$
being the obvious one.
With the notation
$
(q_1,p_1,q_2,p_2) \in (\bold R^3)^4
$
for the elements of $W$,
the momentum mapping $\mu$ is given by the assignment
to $(q_1,p_1,q_2,p_2)$
of $q_1 \wedge p_1 + q_2 \wedge p_2$.
Moreover,
by invariant theory, cf. \cite\weylbook, \cite\locpois,
the ten distinct invariants
$$
q_iq_j,\, q_ip_j,\,p_ip_j,\ 1 \leq i ,j \leq 2,
\tag8.1
$$
among
the scalar products,
together with the four
determinants
$$
|q_1 p_1 q_2|,
\
|q_1 p_1 p_2|,
\
|q_1 q_2 p_2|,
\
|p_1 q_2 p_2|,
\tag8.2
$$
constitute a complete set of invariants
for the
$\roman{SO}(3,\bold R)$-action
on
$W$.
However, for
$(q_1,p_1,q_2,p_2) \in
V= \mu^{-1}(0)$, that is,
when
$
q_1 \wedge p_1
+
q_2 \wedge p_2 = 0,
$
any three of
$(q_1,p_1,q_2,p_2)$
are linearly dependent, that is,
$(q_1,p_1,q_2,p_2)$ lie in a plane in $\bold R^3$, whence
the
four determinants (8.2) vanish
on $V$,
and
the algebra
of invariants
$A[V]^{\roman{SO}(3,\bold R)}$
in the coordinate ring $A[V] = \bold R[W]/I_V$
is in fact generated
by the ten scalar products;
these induce
the quadratic $\roman{SO}(3,\bold R)$-invariant map
$$
\lambda \colon
W @>>>
\roman S^2(\bold R^4),
\quad
\lambda(q_1,p_1,q_2,p_2)
=
\left[
\matrix
q_1 q_1 & q_1 q_2 & q_1 p_1 & q_1 p_2
\\
q_2 q_1 & q_2 q_2 & q_2 p_1 & q_2 p_2
\\
p_1 q_1 & p_1 q_2 & p_1 p_1 & p_1 p_2
\\
p_2 q_1 & p_2 q_2 & p_2 p_1 & p_2 p_2
\endmatrix
\right]
\tag8.3
$$
into the $10$-dimensional real vector space
$\roman S^2(\bold R^4)$
of symmetric 4 by 4 matrices
which, in turn, passes to
an embedding
$$
\widetilde \lambda\colon
W_{\roman{red}}
@>>>
\roman S^2(\bold R^4)
\tag8.4
$$
of
$W_{\roman{red}}$
into
$\roman S^2(\bold R^4)$
as a real semi-algebraic set
$S$;
more details about this semi-algebraic realization will be given below.
We assert at first that
the Zariski tangent space $\roman T_{0}S$ at the origin equals the
whole ambient space, that is, has dimension 10.
In fact, $S$ is a cone since for
$(q_1,p_1,q_2,p_2) \in V$
and $t \in \bold R$,
$$
\widetilde \lambda [t(q_1,p_1,q_2,p_2)]
=
t^2 \widetilde \lambda [q_1,p_1,q_2,p_2] \in S.
$$
Hence for $x \in S$, the half line
$\{t x; t \geq 0\}$ lies in $S$.
Let $v$ be an arbitrary vector in $\bold R^3$ of length one.
Then the vectors
$$
\gathered
(v,0,0,0),\,
(0,v,0,0),\,
(0,0,v,0),\,
(0,0,0,v),\,
(v,v,0,0),
\\
(v,0,v,0),\,
(v,0,0,v),\,
(0,v,v,0),\,
(0,v,0,v),\,
(0,0,v,v)
\endgathered
\tag8.5
$$
all lie in $V$, and
inspection shows that their images
in $S$
under
$\lambda$
are linearly independent
in the ambient vector space
$\roman S^2(\bold R^4)$
and hence constitute a basis.
In fact,
$$
\lambda (v,0,0,0) =
\left[
\matrix
1&0&0&0
\\
0&0&0&0
\\
0&0&0&0
\\
0&0&0&0
\endmatrix
\right],\qquad
\lambda (v,v,0,0) =
\left[
\matrix
1&1&0&0
\\
1&1&0&0
\\
0&0&0&0
\\
0&0&0&0
\endmatrix
\right]
$$
etc.
Consequently the linear span of the cone $S$ equals
the whole ambient space $\roman S^2(\bold R^4)$, and hence the latter
coincides with the
Zariski tangent space $\roman T_{0}S$ at the origin
as asserted.
In particular,
the minimal number of generators of
the algebra $A[V]^{\roman{SO}(3,\bold R)}$
is ten, and this is also the minimal number
of generators of $C^{\infty}(W_{\roman{red}})$
since if fewer generators did suffice
the dimension of the Zariski tangent space would be smaller.
\smallskip
These observations translate to the moduli space $N$ in the following way:
Let
$[A]$  be a point
in
$N_{G}$.
Then the Zariski tangent space
$\roman T_{[A]} N$ has dimension 10 and hence
the minimal number of generators
of
$C^{\infty}(N)$
near $[A]$ or rather that of its germ at $[A]$ is 10.
Moreover, a closer look reveals that
the Zariski tangent space
$\roman T_{[A]} N$
equals that
of
$\roman T_{[A]} N_{(T)}$, with reference to
the induced smooth smooth structure
$C^{\infty}(N_{(T)})$.
In fact,
in the language of constrained systems,
$N_{(T)}$
corresponds to reduced states
where each of the two particles individually
has angular momentum zero,
cf. what is said in our paper
\cite\locpois, and hence
the images
of the ten vectors (8.5)
under
$\lambda$
already lie in the part of $S$
which corresponds to
$N_{(T)}$.
In particular,
the minimal number of generators
of the induced smooth structure
$C^{\infty}(N_{(T)})$
near $[A]$ or rather that of its germ at $[A]$ is still 10.
Finally,
the linear map
$\lambda_A$ from
$\roman H_A^1(\Sigma,\roman{ad}(\xi))$
to
$\roman T_{[A]} N$
is zero
since the derivative of $\lambda$
at the origin
is zero.
Thus
the Zariski tangent space
$\roman T_{[A]} N$
can
in no way be identified with
the cohomology group
$\roman H_A^1(\Sigma,\roman{ad}(\xi))$,
cf. (7.12) above.
It seems also worthwhile
pointing out that,
cf. \cite\locpois,
as a complex variety,
near a point
$[A]$
in
$N_{G}$,
the stratum
$N_{(T)}$
looks like the quadric
$Y^2 = XZ$
in complex 3-space
and hence
at a point
$[A]$
in
$N_{G}$
the complex Zariski
tangent space
of
$N_{(T)}$
has dimension 3.
Thus we see once more that, as a smooth space,
the moduli space $N$ of flat
$\roman {SU}(2)$-connections
for a surface $\Sigma$ of genus 2
looks rather different
from complex projective 3-space
with its standard smooth structure.
\smallskip
More information about the geometry
of $N$ near
a point
$[A]$
in
$N_{G}$
can be obtained in the following way:
The cone $V$ in $W$ may be defined
as the zero locus of the single
homogeneous
real quartic function $\Psi$ on $W$ given by the formula
$$
\Psi(q_1,p_1,q_2,p_2)  =
(q_1 \wedge p_1 + q_2 \wedge p_2)(q_1 \wedge p_1 + q_2 \wedge p_2).
$$
However this function looks like
$$
\Psi(q_1,p_1,q_2,p_2)  =
\left| \matrix q_1 q_1 & q_1 p_1\\
p_1 q_1 & p_1 p_1\endmatrix\right|
+
2\left| \matrix q_1 q_2 & q_1 p_2\\
p_1q_2& p_1 p_2\endmatrix\right|
+\left| \matrix q_2 q_2 & q_2 p_2\\
p_2q_2& p_2 p_2\endmatrix\right|
$$
and hence passes to a quadratic function
$\psi$ on $\roman S^2(\bold R^4)$.
Next we observe that the reduced space
$W_{\roman{red}}$ with respect to
the $\roman {SO}(3,\bold R)$-action coincides
with
the reduced space
with respect to
the action
of the larger group
$\roman {O}(3,\bold R)$
since
the four determinants (8.2) which distinguish between the two
reduced spaces vanish on $V$; this is a special phenomenon
due to the fact that we are considering angular momentum of two particles
in $\bold R^3$.
Now
$W_{\roman{red}}$
appears as the zero locus of the single function
$\psi$ on $W\big /\roman{O}(3,\bold R)$.
However,
by invariant theory,
the ten distinct
inner products (8.1)
constitute a complete set of invariants
for the
$\roman{O}(3,\bold R)$-action
on
$W$
subject to the single defining relation
$$
\left|
\matrix
q_1 q_1 & q_1 q_2 & q_1 p_1 & q_1 p_2
\\
q_2 q_1 & q_2 q_2 & q_2 p_1 & q_2 p_2
\\
p_1 q_1 & p_1 q_2 & p_1 p_1 & p_1 p_2
\\
p_2 q_1 & p_2 q_2 & p_2 p_1 & p_2 p_2
\endmatrix
\right| =0.
\tag8.6
$$
Consequently the affine categorical quotient
$W\big /\big /\roman{O}(3,\bold R)$
amounts to the space of singular symmetric 4 by 4 matrices,
and $W\big /\roman{O}(3,\bold R)$
is realized as its semi-algebraic subset
which consists of non-negative semidefinite matrices.
Thus the reduced space
$W_{\roman{red}}$
and hence the space $N$ near
a point
$[A]$
in
$N_{G}$
appear as the zero locus of the single function
$\psi$ on the subspace of singular non-negative semidefinite matrices.
The determinant and $\psi$
clearly yield
two $\roman{SO}(3,\bold R)$-invariant polynomials
vanishing on $V$,
that is,
elements of the ideal $I_V^{\roman{SO}(3,\bold R)}$
but these two will {\it not\/} generate
$I_V^{\roman{O}(3,\bold R)}$. In fact we can at once write
down
the following six
$\roman{O}(3,\bold R)$-invariant polynomials
which vanish on $V$ and are
quadratic in the generators (8.1) of $A[W]^{\roman{O}(3,\bold R)}$:
$$
(q_1\wedge q_2) \mu,
\quad
(q_1\wedge p_1) \mu,
\quad
(q_1\wedge p_2) \mu,
\quad
(q_2\wedge p_1) \mu,
\quad
(q_2\wedge p_2) \mu,
\quad
(p_1\wedge p_2) \mu
\tag8.7
$$
More explicitly, $(a,b)$
denoting any of the six couples $(q_1,q_2)$ etc.,
we have
$$
\left((a\wedge b) \mu\right)(q_1,p_1,q_2,p_2)
=
\left| \matrix a q_1 & a p_1\\
b q_1 & b p_1\endmatrix\right|
+
\left| \matrix a q_2 & a p_2\\
b q_2 & b p_2\endmatrix\right|.
$$
Moreover, from the six relations
$$
|u_{j_1} u_{j_2} u_{j_3}|
|v_{j_1} v_{j_2} v_{j_3}|
=
\left|
\matrix
u_{j_1} v_{j_1} & u_{j_1} v_{j_2} & u_{j_1} v_{j_3}
\\
u_{j_2} v_{j_1} & u_{j_2} v_{j_2} & u_{j_2} v_{j_3}
\\
u_{j_3} v_{j_1} & u_{j_3} v_{j_2} & u_{j_3} v_{j_3}
\endmatrix
\right|
$$
among the $\roman{SO}(3,\bold R)$-invariants (8.1) and (8.2),
cf. \cite\weylbook,
where $|u_{j_1} u_{j_2} u_{j_3}|$
and
$|v_{j_1} v_{j_2} v_{j_3}|$
refer to any of the four determinants (8.2),
we conclude that
on $V$
all 3 by 3 minors of $\lambda(q_1,p_1,q_2,p_2)$
vanish;
these 3 by 3  minors
yields six additional
${\roman{O}(3,\bold R)}$-invariant polynomials
vanishing on $V$,
of degree three
in the generators
(8.1)
of $A[W]^{\roman{O}(3,\bold R)}$.
In particular, the image $\widetilde \lambda(W_{\roman{red}})$
lies in the subspace of symmetric 4 by 4 matrices
having  rank at most 2.
We conjecture that the six quadratic polynomials
(8.7)
and the six cubic ones
arising from the 3 by 3 minors
constitute a complete set of generators of the ideal
$I_V^{\roman{O}(3,\bold R)}=
I_V^{\roman{SO}(3,\bold R)}$.
\smallskip
The methods of
{\smc Lerman-Montgomery-Sjamaar}~\cite\lermonsj\
yield a geometric description
of $W_{\roman{red}}$,
viewed as a
subspace
of that of symmetric 4 by 4 matrices:
Let $J$ be the symplectic operator on $\bold R^4$:
$J^2 = -1,\, J\,{}^{\roman{t}}\negthinspace J = \roman{Id},\,
\sigma(v,w) = v Jw$.
The
assignment
$S \mapsto JS$
identifies
$\roman S^2(\bold R^4)$
with the
Lie algebra $\roman {sp}(2,\bold R)$,
and
a result in \cite\lermonsj\
implies that
$\widetilde \lambda$ identifies
$W_{\roman{red}}$
with the closure of the nilpotent orbit
in $\roman {sp}(2,\bold R)$
which corresponds to
positive symmetric 4 by 4 matrices
of rank at most 2 having kernel a coisotropic subspace.
The Lie algebra
$\roman {sp}(2,\bold R)$
has rank two --- in fact it is the split real from of
$\roman C_2$ which coincides with $\roman B_2$, though ---
and its algebra of
$\roman {Sp}(2,\bold R)$-invariants
under the adjoint representation
is a polynomial algebra, generated by the Killing form
and the determinant.
Hence the nilvariety $\roman{Nil}(\roman {sp}(2,\bold R))$
is of real dimension 8; it
consists of singular matrices in
$\roman {sp}(2,\bold R)$
having vanishing Killing form, and its subspace
$\roman{Nil}^+(\roman {sp}(2,\bold R))$
of non-negative semidefinite
matrices
is a union
$\bold n_0 \cup \bold n_1
\cup \bold n_2
\cup \bold n_3$
of four
nilpotent adjoint orbits,
$\bold n_j$ being
the subspace of non-negative semidefinite
rank $j$ matrices.
The reduced space
$W_{\roman{red}}$
now appears as the
union
$\bold n_0 \cup \bold n_1
\cup \bold n_2$.
It may be described as a zero locus
in $\roman{Nil}^+(\roman {sp}(2,\bold R))$
in various ways, that is,
\newline\noindent
--- of the function $\psi$ or what corresponds to it,
restricted to
$\roman{Nil}^+(\roman {sp}(2,\bold R))$,
\newline\noindent
--- of the  functions (8.7) or what corresponds to them,
restricted to
$\roman{Nil}^+(\roman {sp}(2,\bold R))$;
in fact, the two functions $(q_1 \wedge p_1)\mu$ and
$(q_2 \wedge p_2)\mu$ already suffice;
\newline\noindent
--- of the six 3 by 3 minors,
restricted to
$\roman{Nil}^+(\roman {sp}(2,\bold R))$.
\newline\noindent
Somewhat amazingly,
since, as a space, $W_{\roman{red}}$ is smooth in the ordinary sense,
in fact a copy of
real affine 6-dimensional space,
the
union
$\bold n_0 \cup \bold n_1
\cup \bold n_2$
is just a real affine 6-dimensional space.

\bigskip
\centerline{\smc References}
\medskip
\widestnumber\key{999}
\ref \no  \abramars
\by R. Abraham and J. E. Marsden
\book Foundations of Mechanics
%\bookinfo Graduate Texts in Mathematics, No. 60
\publ Benjamin/Cum-\linebreak
mings Publishing Company
%\publaddr
\yr 1978
\endref

\ref \no  \armgojen
\by J. M. Arms, M. J. Gotay, and G. Jennings
\paper  Geometric and algebraic reduction for singular
momentum mappings
\jour Advances in Mathematics
\vol 79
\yr 1990
\pages  43--103
\endref

\ref \no  \armamonc
\by J. M. Arms, J. E. Marsden, and V. Moncrief
\paper  Symmetry and bifurcation of moment mappings
\jour Comm. Math. Phys.
\vol 78
\yr 1981
\pages  455--478
\endref

\ref \no  \atibottw
\by M. Atiyah and R. Bott
\paper The Yang-Mills equations over Riemann surfaces
\jour Phil. Trans. R. Soc. London  A
\vol 308
\yr 1982
\pages  523--615
\endref

\ref \no \bocosroy
\by J. Bochnak, M. Coste, and M.-F. Roy
\book G\'eom\'etrie
alg\'ebrique r\'eelle
\bookinfo Ergebnisse der Mathematik und ihrer Grenzgebiete, 3. Folge, Band 12.
A series of modern surveys in Mathematics
\vol 12
\yr 1987
\publ Springer
\publaddr Berlin-Heidelberg-New York-Tokyo
\endref

\ref \no  \ebinone
\by D. Ebin
\paper The manifold of Riemannian metrics
\jour Proceedings of Symposia in Pure Mathematics
\vol 15
\yr 1970
\pages 11--40
\publ Amer. Math. Soc.
\publaddr Providence R. I
\endref

\ref \no  \eellsone
\by J. Eells
\paper A setting for global analysis
\jour Bull. Amer. Math. Soc.
\vol 72
\yr 1966
\pages 751--807
\endref

\ref \no  \goldmone
\by W. M. Goldman
\paper The symplectic nature of the fundamental group of surfaces
\jour Advances
\vol 54
\yr 1984
\pages 200--225
\endref

\ref \no \howeone
\by R. Howe
\paper Remarks on classical invariant theory
\jour  Trans. Amer. Math. Soc.
\vol 313
\yr 1989
\pages  539--570
\endref

\ref \no \singula
\by J. Huebschmann
\paper The singularities of Yang-Mills connections
for bundles
on a surface. I. The local model
\paperinfo Math. Z. (to appear), dg-ga/9411006
%\jour \vol \yr \pages
\endref
\ref \no  \singulat
\by J. Huebschmann
\paper The singularities of Yang-Mills connections
for bundles on a surface. II. The stratification
\paperinfo Math. Z. (to appear), dg-ga/9411007
%\jour \vol \yr \pages
\endref
\ref \no \topology
\by J. Huebschmann
\paper
Holonomies of central Yang-Mills connections
for bundles  on a surface
with disconnected structure group
\jour Math. Proc. Camb. Phil. Soc.
\yr 1994
\vol 116
\pages 375--384
\endref

\ref \no \direct
\by J. Huebschmann
\paper
Poisson
structures on certain
moduli spaces
for bundles on a surface
\paperinfo Annales de l'Institut Fourier (to appear)
\endref

\ref \no \singulth
\by J. Huebschmann
\paper The singularities of Yang-Mills connections
for bundles on a surface. III. The identification of the strata
\paperinfo in preparation
%\jour \vol \yr \pages
\endref

\ref \no \locpois
\by J. Huebschmann
\paper Poisson geometry of
flat connections for {\rm SU(2)}-bundles on surfaces
\paperinfo Math. Z. (to appear), hep-th/9312113
%\jour \vol \yr \pages
\endref

\ref \no \modusym
\by J. Huebschmann
\paper Symplectic and Poisson structures of certain moduli spaces
\paperinfo Preprint 1993, hep-th/9312112
%\jour \vol \yr \pages
\endref

\ref \no \huebjeff
\by J. Huebschmann and L. Jeffrey
\paper Group cohomology construction
of symplectic forms on certain moduli spaces
\jour Int. Math. Research Notices
\vol 6
\yr 1994
\pages 245--249
\endref

\ref \no \kobanomi
\by S. Kobayashi and K. Nomizu
\book Foundations of differential geometry, I (1963), II (1969)
\bookinfo Interscience Tracts in Pure and Applied Mathematics, No. 15
\publ Interscience Publ.
\publaddr New York-London-Sydney
%\yr 1963
\endref

\ref \no \kostafou
\by B. Kostant
\paper Lie group representations on polynomial rings
\jour Amer. J. of Math.
\vol 85
\yr 1963
\pages 327--404
\endref

\ref \no \lermonsj
\by E. Lerman, R. Montgomery and R. Sjamaar
\paper Examples of singular reduction
\book Symplectic Geometry
\bookinfo Warwick, 1990, ed. D. A. Salamon, London Math. Soc. Lecture Note
Series
\vol 192
\yr 1993
\pages  127--155
\publ Cambridge University Press
\publaddr Cambridge, UK
\endref

\ref \no \marswein
\by J. Marsden and A. Weinstein
\paper Reduction of symplectic manifolds with symmetries
\jour Rep. on Math. Phys.
\vol 5
\yr 1974
\pages 121--130
\endref

\ref \no \mathone
\by J. Mather
\paper Differentiable invariants
\jour Topology
\vol 16
\yr 1977
\pages 145--155
\endref

\ref \no \mittvial
\by P. K. Mitter and C. M. Viallet
\paper On the bundle of connections and the gauge orbit manifold in
Yang-Mills theory
\jour Comm. in Math. Phys.
\vol 79
\yr 1981
\pages 457--472
\endref

\ref \no \narashed
\by M. S. Narasimhan and C. S. Seshadri
\paper Stable and unitary vector bundles on a compact Riemann surface
\jour Ann. of Math.
\vol 82
\yr 1965
\pages  540--567
\endref

\ref \no \naramntw
\by M. S. Narasimhan and S. Ramanan
\paper Moduli of vector bundles on a compact Riemann surface
\jour Ann. of Math.
\vol 89
\yr 1969
\pages  19--51
\endref

\ref \no \nararama
\by M. S. Narasimhan and T. R. Ramadas
\paper Geometry of $\roman{SU}(2)$ gauge fields
\jour Comm. in Math. Phys.
\vol 67
\yr 1979
\pages  121--136
\endref

\ref \no \raghuboo
\by M. S. Raghunatan
\book Discrete subgroups of Lie groups
%\bookinfo Ergebnisse dr Mathematik
\publ Springer
\publaddr Berlin-\linebreak
Heidelberg-New York
\yr 1972
\endref

\ref \no \gwschwar
\by G. W. Schwarz
\paper Smooth functions invariant under the action of
a compact Lie group
\jour Topology
\vol 14
\yr 1975
\pages 63--68
\endref

\ref \no \gwschwat
\by G. W. Schwarz
\paper The topology of algebraic quotients
\paperinfo In: Topological methods in algebraic
transformation groups
\jour Progress in Math.
\vol 80
\yr 1989
\pages 135--152
\publ Birkh\"auser
\endref

\ref \no \sjamlerm
\by R. Sjamaar and E. Lerman
\paper Stratified symplectic spaces and reduction
\jour Ann. of Math.
\vol 134
\yr 1991
\pages 375--422
\endref

\ref \no  \weilone
\by A. Weil
\paper  On discrete subgroups of Lie groups
\jour
(a) Ann. of Math.
\vol 72
\yr 1960
\pages  369--384
\moreref
(b) Ann. of Math.
\vol 75
\yr 1962
\pages  578--602
\endref

\ref \no  \weiltwo
\by A. Weil
\paper  Remarks on the cohomology of groups
\jour
Ann. of Math.
\vol 80
\yr 1964
\pages  149--157
\endref

\ref \no  \weylbook
\by H. Weyl
\book The classical groups
%\bookinfo
\publ Princeton University  Press
\publaddr Princeton NJ
\yr 1946
\endref

\ref \no  \whitnone
\by H. Whitney
\paper Analytic extensions of differentiable functions defined
on closed sets
\jour Trans. Amer. Math. Soc.
\vol 36
\yr 1934
\pages  63--89
\endref

\ref \no  \whitnfou
\by H. Whitney
\paper Elementary structure of real algebraic varieties
\jour  Ann. of Math.
\vol 66
\yr 1957
\pages  545--556
\endref

\enddocument